\documentstyle[11pt]{article}
\begin{document}
\baselineskip 16pt
\begin{center}
{\large \bf ABS ALGORITHMS FOR LINEAR SYSTEMS
		     AND OPTIMIZATION:\\
A REVIEW AND A BIBLIOGRAPHY}
\vskip8mm
Emilio Spedicato$^*$, \ Elena Bodon$^*$, \   Antonino Del Popolo$^*$ \
and \ Zunquan Xia$^{+}$
\vskip5mm
* -- Department of Mathematics, University of Bergamo \\
  + -- Department of Applied Mathematics, Dalian University of Technology
\end{center}
\vskip8mm

\noindent
{\bf 1. The scaled ABS class}.
\vskip3mm

\noindent ABS methods have been introduced by Abaffy, Broyden and
Spedicato (1984), in a paper where the solution of linear algebraic
equations was considered via the so called
 {\it ABS unscaled or basic class}. That class was later generalized to the so
 called 
  {\it scaled ABS class} and then applied also to the solution of
linear least squares, nonlinear algebraic equations and optimization problems.
For a full presentation of the main results in the ABS field see the monographs
of Abaffy and Spedicato (1989) and of  Zhang, Xia and Feng (1999).
There are presently almost 400
 papers dealing with ABS methods, see for a partial list the bibliography
 in the Appendix 2. In  this review we restrict our attention to
 ABS methods for linear determined or underdetermined algebraic systems
 and their application to some optimization problems. We begin by
 giving our notation and the steps of the scaled ABS class.
	 
 Consider the following general linear system (determined or underdetermined),
where the rank of $A \in R^{m,n}$ is arbitrary

$$
Ax = b \ \ \ \ \  x \in {R}^n ,\ \  b \in {R}^m , \ \ m \leq n \eqno {(1)}
$$
or
$$
a_i^Tx - b_i = 0, \ \ \  i = 1,\ldots,m  \eqno {(2)}
$$
where
$$
A = \left
    [
    \begin{array} {c}
    a_1^T\\
    \cdots\\
    a_m^T
    \end{array}
    \right ] \eqno {(3)}
$$
Then the scaled ABS class is defined by the following procedure:

\begin{enumerate}
\item Let $x_1 \in {R}^n$ be arbitrary,\ let $H_1 \in {R}^{n,n}$ 
be nonsingular arbitrary, let $v_1 \in R^m$ be arbitrary nonzero, set
$i = 1$.
\item Compute the residual vector $r_i = Ax_i-b$. 
If $r_i = 0$ stop ($x_i$
solves the system), otherwise compute
$s_i = H_iA^Tv_i$. If $s_i \neq 0$, go to (3).  If $s_i = 0$
and $\tau = v_i^Tr_i = 0$, set 
$x_{i+1} = x_{i} ,\  H_{i+1} = H_{i}$
and go to  (6). Otherwise  stop, the system has no solution.
\item Compute the search vector $p_i$ by
$$
p_i = H_i^Tz_i \eqno{(4)}
$$
with $z_i \in  {R}^n$ arbitrary save that
$$
v_i^TAH_i^Tz_i \neq 0 \eqno{(5)}
$$
\item Update the estimate of the solution by
$$
x_{i+1} = x_{i} - \alpha_ip_i  \eqno{(6)}
$$
where the stepsize $\alpha_i$  is given by
$$
\alpha_i = v_i^Tr_i/v_i^TAp_i. \eqno{(7)}
$$
\item Update the matrix $H_i$ by
$$
H_{i+1} = H_i - H_iA^Tv_iw_i^TH_i/w_i^TH_iA^Tv_i \eqno{(8)}
$$
where $w_i \in {R}^n$ is arbitrary save that
$$
w_i^TH_iA^Tv_i \neq 0.   \eqno{(9)}
$$
\item If $i = m$, then stop  ($x_{m+1}$ solves the system), otherwise
give $v_{i+1} \in R^m$ arbitrary linearly
independent from $v_1, \ldots, v_i$, increase
$i$ by one and go to (2).
\end{enumerate}
Matrices $H_i$ are generalizations of oblique projectors. They were first
introduced in literature apparently by  Wedderburn (1934) and
have been used in several papers of Egervary  that
in fact introduced ABS methods ante litteram. They have been denoted as
 {\it  Abaffian matrices} at the First International Conference on ABS Methods
 (Luoyang, China, 1991) and this name will be used here.

The above recursions define a class of algorithms. A particular method is
obtained by a choice of the parameters $H_1, \ v_i, \ ,z_i, \  w_i$.
The {\it basic ABS class} is obtained by setting
 $v_i = e_i, \ e_i$ being the
$i$-th unit vector  in $R^m$. Parameters $w_i, z_i, H_1$ 
were introduced respectively by 
  Abaffy, Broyden and Spedicato, from whose initials the class name
is obtained. One can show that the scaled ABS class is a complete
realization of the so called 
 Petrov-Galerkin principle for solving a linear system,
where the iteration has the form $x_{i+1} = x_i - \alpha_ip_i$
with $\alpha_i,\ p_i$ such that
$r_{i+1}^{\ T}v_j = 0, \ j=1, \ldots, i$, where vectors $v_j$ are
arbitrary linearly independent. It appears that all methods proposed
 in the literature for solving a linear system in a finite number
of steps are special cases of the ABS class. Indeed even the
  Quasi-Newton methods of the  Broyden's class for nonlinear systems, which
are known to have finite termination under mild conditions  
in no more that $2n$ steps on linear systems, have been shown recently
to be imbeddable in the ABS class, see Adib, Mahdavi Amiri and Spedicato (1999), 
in the sense that the approximation to the solution of index
 $2i-1$ generated by Broyden method is identical with the approximation
of index $i$ of a certain ABS algorithm, when both algorithms start
from a same initial point.

The recursions of the scaled ABS class can be derived in different ways.
For instance, let us consider the basic ABS class, which has the property that
the approximation $x_{i+1}$ solves the first $i$ equations of the system
$Ax=b$.   This property is obtained by a suitable choice of the stepsize
and the search vector in the update $x_{i+1} = x_i - \alpha_ip_i$. Choosing
$\alpha_i$ by formula (7) makes the $i-$th equation $a_i^Tx = b_i$ to be
satisfied by $x_{i+1}$. Choosing $p_i$ by (4) with $H_i$ updated by (8)
makes the previous $i-1$ equations to remain satisfied in the new point,
since a simple induction argument shows that this is true if $H_ia_j = 0$
for $j < i$, a property that characterizes update (8). Conditions (5) and (9)
guarantee that the process is well-defined (no division by zero).
Condition (5) has a precise geometrical justification, implying that the
$i-$th search vector intersects the linear variety containing all
solutions of the $i-$th equation at a finite point (the vector $a_i$ being
orthogonal to such a variety). The recursions of the scaled ABS class are
obtained from those of the basic ABS class by applying them on the
system $V^TAx = V^Tb$, where $V = (v_1, \ldots, v_m)$, which is equivalent in
the sense of possessing the same set of solutions, if $V$ is nonsingular.

Referring for proofs to the monograph of Abaffy and 
Spedicato (1989), we recall now some properties
of the scaled ABS class. For simplicity we assume that $rank(A)=m$.
\begin{itemize}
\item Define  $V_i = (v_1,\ldots,v_i)$,
$W_i = (w_1,\ldots,w_i)$. Then 
$H_{i+1}A^TV_i = 0,\  H_{i+1}^TW_i = 0$,  hence
vectors  $A^Tv_j, w_j, \ j=1,\ldots,i$ span respectively
$Null(H_{i+1})$ and $Null(H_{i+1})^T$.
\item Vectors $H_iA^Tv_i, \  H_i^Tw_i $ are zero iff
 $a_i, \ w_i$  respectively linearly depend on
$a_1,\ldots,a_{i-1}$, $w_1, \ldots, w_{i-1}$.
\item Define  $P_i = (p_1,\ldots,p_i) $ and $A_i = (a_1,\ldots,a_i)$.
Then we have the implicit factorization
$V_i^TA_i^TP_i = L_i$, with $L_i$  
nonsingular lower triangular. Hence for
 $m=n$ we get a semiexplicit factorization of the inverse of
 $A$, i.e. with $P = P_n, \ V = V_n, \ L = L_n$
$$
A^{-1} = PL^{-1}V^T.  \eqno{(10)}
$$
For several choices of  $V$, matrix  $L$ is diagonal, hence
formula (10) gives an explicit factorization of the inverse, a bonus freely
provided by the ABS procedure. One can moreover show that all possible
factorizations of a matrix can be obtained by suitable choices of the ABS 
parameters, another completeness result.

\item Define $S_i$ and $R_i$ as $ S_i = (s_1,\ldots,s_i), \
R_i = (r_1,\ldots,r_i) $, where
$s_i = H_iA^Tv_i, \ r_i = H_i^Tw_i $.
Then
$H_{i+1} = H_1 - S_iR_i^T$
and vectors $s_i$, $r_i$  and $p_i$ can be built by a generalized
 Gram-Schmidt iteration using $2i$ vectors. This alternative formulation is computationally
convenient in term of storage when $m < n/2$. There is also an alternative
 formulation using $n-i$ vectors.

\item One can write the  Abaffian as follows, in term of the parameter
matrices
$$
H_{i+1} = H_1 - H_1A^TV_i(W_i^TH_1A^TV_i)^{-1}W_i^TH_1. \eqno{(11)}
$$
Setting $V = V_m, \ W = W_m$, one can prove that 
$H_1, \ 
V, \ W$ are admissible (i.e. condition (9) is satisfied) iff matrix
 $Q = V^TAH_1^TW$ is strongly nonsingular (i.e. $Q$ is LU
factorizable).

\item If $m < n$ the general solution of $Ax=b$ can be written as follow,
with  $q \in R^n$ arbitrary
$$
x = x_{m+1} + H^T_{m+1}q       \eqno{(12)}
$$
\end{itemize}

It was observed that the {\it scaled ABS class} can be obtained by
applying the unscaled ABS class to the scaled system
$V^TAx = V^Tb$. This scaled system can be viewed also as a left
preconditioned system. Hence the scaled ABS class includes also all
 possible left preconditioners, which are defined in the ABS context 
implicitly and dynamically.
\vskip8mm
\noindent {\bf 2. Some subclasses and methods}   
\vskip5mm
\noindent Here we consider some subclasses and special methods in the scaled ABS class.
\begin {description}
\item[(a)] The {\it conjugate direction subclass}. \ \
This subclass is given by the choice $v_i=p_i$ and is well defined under the
sufficient but not necessary condition that
 $A$ is symmetric positive definite. All methods in this subclass have
 search directions that are $A$-conjugate and conversely any method that
generates $A$-conjugate directions is a member of the subclass. If $x_1=0$ then
$x_{i+1}$ is the vector of least energy norm (i.e. $A$-weighted Euclidean
norm) in the linear variety
$x_1 + Span(p_1,\ldots,p_i)$ and moreover the solution is approached monotonically
in such a norm. The subclass contains the ABS reformulation of the
Choleski algorithm and of the methods of
 Hestenes-Stiefel and of Lanczos. 

\item[(b)] The {\it orthogonally scaled subclass}. \ \
This subclass is obtained by the choice $v_i=Ap_i$ and is well defined if
$A$ has full column rank; it is also well defined even if
 $m > n$, in which case the iteration can be performed for $n$ steps and
 will stop in a least squares solution. The subclass generates scaling
 vectors that are orthogonal and search directions that are $A^TA$-conjugate.
If $x_1=0$ then
 $x_{i+1}$ minimizes the Euclidean norm of the residual in the linear
variety $x_1 + Span(p_1, \ldots,p_i)$, implying that the solution is
monotonically approached from below in the residual norm. This subclass
contains the ABS formulation of the QR and GMRES methods.

\item[(c)] The {\it optimally stable subclass}. \ \ This subclass is
obtained by the choice $v_i=A^{-T}p_i$, the inverse disappearing in the
actual recursion formulas. The class is well defined without conditions
on the matrix (assumed to be square nonsingular, a condition that can be
relaxed). The search vectors in the subclass are orthogonal.
If $x_1=0$ then
 $x_{i+1}$ is the vector of least  Euclidean norm  that solves the first
$i$ equations, implying that the solution is
monotonically approached from below in the Euclidean norm. The 
name of the subclass derives from the fact that a condition is satisfied
which guarantees the least propagation of a single error introduced in
the sequence $x_i$. This subclass contains the Huang algorithm, see
the following paragraph, and the ABS version of the Craig algorithm.
\item[(d)] The {\it Huang method.} \ \ This method belongs to both the
basic ABS class and the optimally scaled subclass. It can be obtained by
the parameter choices   $H_1 = I$, $z_i =
w_i = a_i, \ v_i = e_i$. A mathematically equivalent, in the sense of generating
the same iterates in exact arithmetic, but numerically more stable form is
the {\it modified Huang method}, based upon formulas
 $p_i = H_i(H_ia_i)$ and
$H_{i+1} = H_i - p_ip_i^T/p_i^Tp_i$. The Huang method generates orthogonal
search vectors. If $x_1=0$ then  $x_{i+1}$ is the solution of least
 Euclidean norm of the first $i$ equations and the solution $x^+$ of
least  Euclidean norm of the system is approached monotonically and from below
from sequence  $x_i$. Zhang (1995a) has shown that using the active set
strategy of  Goldfarb and Idnani (1983) the Huang algorithm can solve a
system of general linear inequalities where in a finite number of steps
it either establishes that no feasible point exists or it finds the solution
of least Euclidean norm. If one takes as initial matrix a symmetric positive
definite matrix $B^{-1}$, then all above properties remain valid with respect
to the $B$-weighted Euclidean norm.
 
Extensive numerical experiments, see Spedicato, Bodon and Luksan (2000) have
shown that a certain version of the modified Huang algorithm is very
accurate, being comparable on extremely ill conditioned problems with
the LAPACK implementations of the rank revealing QR factorization or the
SVD factorization. It 
provides usually the same estimate of the numerical rank
as these  LAPACK codes. Its computational times are  up to
over one hundred times lower on problems with 2000 variables.

\item[(e)] The {\it  implicit LU method. } This method is obtained by the
choices $H_1 = I,\  z_i = w_i = v_i = e_i$, which are well defined iff
 $A$ is strongly nonsingular (otherwise column pivoting is needed, or
 row pivoting if $m = n$). The parameter choices imply the following
structure for the Abaffian, with
$K_i \in R^{n-i,i}$
$$H_{i+1} = \left
	  [
	  \begin{array}{cc}
	  0      & 0 \\
	  \cdots & \cdots \\
	  0      & 0 \\
	  K_i    & I_{n-i}
	  \end{array}
	  \right ]. \eqno{(13)}
$$
From (13) we get that the last $n-i$ components of $p_i$  are zero,
while the $i$-th component equals one, hence, if $x_1=0$, 
the sequence $x_i$ consists of special basic vectors. Moreover
$P_i$ is upper triangular with units on diagonal, hence the associated
implicit factorization is of the LU type, justifying the name. If $m=n$,
then $n^3/3$ multiplications are required, plus lower order terms. Main 
storage is for $K_i$, which is at most
 $n^2/4$. Hence the algorithm can be considered more efficient than
 classical Gaussian elimination or LU factorization, since it requires
 the same number of multiplications but has a lower storage need (in the
case that $A$ cannot be overwritten).

\item[(f)] The {\it implicit LX method } is given by choices
$H_1 = I,v_i = e_i,\ z_i = w_i = e_{k_i}$,
with $k_i$, \ $ 1 \leq k_i \leq n$ such that
$$
e_{k_i}^TH_ia_i \neq 0. \eqno{(14)}
$$ 
If   $A$ has full rank, then there exists at least one
index $k_i$ such that (14) is satisfied. For stability we usually choose
  $k_i$ to maximize
 $  \mid e_{k_i}^TH_ia_i \mid$.

 Let  
$N = (1,2,\ldots,n)$, let $B_i =(k_1,\ldots,k_i)$ and let
 $N_i =N\backslash B_i$. Then we have:
\begin{itemize}
\item  $k_i \in N_{i-1}$
\item  $H_{i+1}^Te_k=0$ for $k \in{{ B_i}}$ 
\item Vector $p_i$  has $n-i$ zero components; the 
$k_i$-th component equals one
\item If $x_1 = 0$, then $x_{i+1}$ is a basic type solution of first
$i$ equations
\item Columns of $H_{i+1}$  of index
 $k \in$ {\bf $N_i$} are the unit vectors $e_k$; 
columns of index $k \in$ {\bf $B_i$} 
have zero components in
 $j$-th position, with $j \in$ {\bf $B_i$}. 
\item At  $i$-th step  $i(n-i)$ multiplications are required to compute
$H_ia_i$ and $i(n-i)$ to update the non trivial part of $H_i$.
Hence the operation cost is the same as for the implicit 
LU algorithm (i.e. $n^3/3$), but pivoting is not required.
\item Storage is again the same as for the implicit LU algorithm, i.e. at
most $n^2/4$. 
\end{itemize}
It follows from the above that the implicit LX algorithm can be
considered as {\it the most efficient general solver} for a linear system not of
the Strassen type, since it has the same operation cost of Gaussian elimination,
requires less storage (when $A$ is not overwritten) and requires no condition
on the matrix. It should also be remarked that
numerical experiments by  Mirnia (1996) 
have shown that it provides a very accurate estimate of the solution,
usually better than by the implicit LU methods or the classical LU
factorization as implemented by MATLAB or LAPACK.
\end{description}
\vskip8mm
{\bf 3. Other ABS algorithms and numerical experiments}  
\vskip5mm
\noindent For most algorithms proposed in the literature their ABS reformulation
has been obtained, leading, in view of the several existing alternative
formulations of the ABS linear algebra, to new formulations in many cases
with better computation cost, storage requirement, numerical stability
and possibility of vectorization and parallelization. The reprojection
technique used in the modified Huang algorithm with success, based on the
identities 
$H_iq = H_i(H_iq), \ H_i^Tq = H_i^T(H_i^Tq)$, valid for any vector $q$
if $H_1 = I$, has been used to improve substantially the stability of 
the methods of Hestenes-Stiefel and of Craig. Moreover the classical
 iterative refinement can be reformulated in the ABS context using
the previous search vectors. Experiments by
 Mirnia (1996) have shown that such refinement coupled with residual
computation in double precision provides usually a fully accurate single
precision solution in just one or two refinement steps.

For structured systems one can  develop special ABS methods exploiting
the problem structure; often the Abaffian matrix reflects the structure
 of $A$. For instance if
 $A$ is banded, so is the Abaffian generated by implicit LU, implicit QR and
 Huang methods, albeit with a larger band width. If
 $A$ is symmetric positive definite with a Nested Dissection structure, so
is the Abaffian generated by the implcit Choleski algorithm, whose storage
turns out to be, for $n$ going to infinity, infinitesimal with respect
to the storage required for the corresponding special implementation of the
classical  Choleski method. For matrices having the
KT \  (Kuhn-Tucker) structure, large classes of ABS methods have been
derived, many of them having up to one order less storage or number
of operations than the classical methods based upon LU or QR factorization
or the Aasen process, see
 Spedicato, Chen and Bodon (1996). Moreover testing of some of such ABS
 methods has also shown that they are  usually more accurate, see
 Spedicato, Bodon and Luksan (2000).
 
 For matrices with a general sparsity pattern little is  known about
the fill-in structure of the Abaffian. We have however found that by
clever use of BLAS4 routines the performance of the full version is greatly
improved and becomes competitive with the performance of codes like MA28
minimizing the fill-in in LU factors, at least for $n$ not too large.

There are also block formulations of the ABS procedure, which do not lead to
accuracy deterioration and provide reduction in number of operations and,
in particular, faster implementations on vector or parallel computers.
The special case of the two-block has been investigated by Adib and 
Mahdavi-Amiri (2000), with good numerical results.

Finally ABS methods can be used to generate large classes of iterative methods
for linear systems via two approaches. The first is based upon restarting 
the iteration after 
$k < m$ steps and using the vector formulation instead of the full Abaffians,
so that the storage is only of order $2k$ vectors.   The second approach
 is based upon the truncation idea. The vector formulation contains
 summations, at the $i-$th step from 1 to $i$. Thus we only keep $p$
 terms in such summations, typically the last $p$ terms, so that storage
 is of order $p$ vectors. Global convergence at a $Q$-linear rate has been
 established for large subclasses of the scaled ABS class, including the
 conjugate direction, the orthogonally scaled and the optimally stable subclasses,
 under mild conditions on the parameter choices, following the technique of
  Dennis and Turner (1987). The obtained class of iterative methods
  includes such classical methods as the
  Krylov space methods, the Gauss-Seidel, the
De La Garza and Kackmartz methods, see  Spedicato and Li (2000).

A project named ABSPACK is now under development aiming at producing
documented software for solving linear and nonlinear algebraic equations
and some optimization problems using ABS methods. In the course of 1999
algorithms have been implemented for solving full linear underdetermined, determined
and overdetermined (in the least squares sense) and KT 
systems. Comparisons have been
made with other commercial packages, including LINPACK, LAPACK, MATLAB, UFO,
showing that the produced codes are competitive in terms of both accuracy
and timing (better results in timing are expected with the use of the
block formulations, currently under development). In the Appendix some
results are given, taken from Spedicato, Bodon and Luksan (2000).
\vskip8mm
\noindent {\bf 4. ABS solution of linear Diophantine equations}
\vskip5mm
\noindent The ABS algorithm has been recently applied to solve linear Diophantine
 systems, where it  determines
if the  system has an
integer solution,  computes a particular solution and provides a
representation of  all integer solutions. The ABS algorithm 
is a generalization of
  a method proposed by Egervary (1955) for the 
 particular case of a homogeneous system.

Let ${Z}$ be the set of all integers and consider the Diophantine linear
system of equations
$$
Ax=b, \quad  \quad x\in {Z}^n,             \
A\in {Z}^{m\times n}, \ b\in {Z}^m, \ m\leq n.\eqno{(15)}
$$

We  recall some results from number theory.
Let $a$ and $b$ be integers. If there is integer $\gamma$ 
so that $b= \gamma a$, then
we say that $a$ divides $b$ and write $a|b$, otherwise 
we write $a \!\!\not|\ b$.
If $a_1,\dots,a_n$ are integers, not all being zero, then the greatest common
divisor ($gcd$) of these numbers is the greatest positive integer $\delta$ which
divides all
$a_i$, $i=1,\dots,n$ and we write $ \delta =gcd(a_1,\dots,a_n)$. We
note that
$\delta \geq 1$ and that 
$\delta$ can be written as an integer linear combination of
the $a_i$, i.e. $\delta = z^Ta$ for some $z\in R^n$.
One can show that $\delta$ is the least positive integer for which
the equation $a_1x_1+\dots+a_nx_n= \delta $ has an integer solution. Now
$\delta$ plays a main role in the following
\vskip2mm
	 
\noindent
{{\bf  Fundamental Theorem of the  
Linear Diophantine Equation}}
\vskip2mm
{\it
\noindent Let $a_1,\dots,a_n$ and $b$ be integer numbers. Then the
 Diophantine linear equation $a_1x_1+\dots+a_nx_n=b$
has  integer solutions if and only
if $gcd(a_1,\dots,a_n) |\ b$. 
In such a case if $n>1$ then
there is an infinite number of integer solutions. }

\vskip2mm

\noindent In order to find the general integer solution of the Diophantine
equation $a_1x_1+\dots+a_nx_n=b$, the main step is to solve
$a_1x_1+\dots+a_nx_n= \delta$,
where $\delta=gcd(a_1,\dots,a_n)$, for a special integer solution. 
There exist several
algorithms for this problem.
The basic step is the computation of $\delta$ and $z$,
often done using  the  algorithm of
Rosser (1941), which  
 avoids a too rapid growth of the intermediate integers, and which
 terminates in polynomial time, as shown
by Schrijver (1986). The scaled ABS algorithm can be applied to Diophantine  equations via a special
choice of its parameters, originating from the following considerations
and Theorems.

Suppose $x_i$ is an integer vector (we begin with $x_1$ an integer vector).
Since $x_{i+1}=x_i-\alpha_ip_i$, then $x_{i+1}$ is integer if
$\alpha_i$ and $p_i$ are integers. If $v_i^TAp_i|(v_i^Tr_i)$
 then $\alpha_i$ is an integer. If $H_i$
and $z_i$ are respectively an integer
matrix and an integer vector, 
then $p_i=H_i^Tz_i$ is also an integer vector. Assume
$H_i$ is an integer matrix. From (6), if
$v_i^TAH_i^Tw_i$ divides all the components of $H_iA^Tv_iw_i^TH_i$, then $H_{i+1}$
is an integer matrix. 

Conditions for the existence and determination of the 
general solution of the Diophantine system are given in the following
theorems,  generalizing the Fundamental Theorem, see 
Esmaeili, Mahdavi-Amiri and Spedicato
 (1999).

\vskip2mm

\noindent {\bf Theorem 1} 
\vskip2mm
\noindent {\it Let $A$ be full rank and suppose that the Diophantine
system (15) is integerly solvable. Consider the  Abaffians generated
by the scaled ABS algorithm with the  parameter choices:
  $H_1$ is unimodular (i.e. both  $H_1$ and $H_1^{-1}$ are integer
  matrices);
  for $i = 1, \ldots, m$, \  $w_i$ is such that
$w_i^TH_iA^Tv_i = \delta_i,  \ \ \delta_i = gcd(H_iA^Tv_i)$.
Then the following properties are true:
\begin{description}
\item[(a)] the  Abaffians generated by the algorithm are
well-defined and are integer matrices
\item[(b)] if $y$ is a special integer solution of the first $i$
equations, then any integer solution $x$ of such equations can be
written as $x = y + H_{i+1}^Tq$ for some integer vector
$q$.
\end{description}    }
\vskip2mm
  
\noindent {\bf Theorem 2} 
\vskip2mm
{\it 
\noindent Let $A$ be full rank and consider the sequence of
matrices $H_i$ generated by the scaled ABS algorithm with parameter
choices as in Theorem 1. Let  $x_1$ in the scaled
ABS algorithm be an arbitrary integer vector and let $z_i$ be  such that
$z_i^TH_iA^Tv_i = gcd(H_iA^Tv_i)$.
Then system (15)     has integer solutions iff 
$gcd(H_iA^Tv_i)$ divides $v_i^Tr_i $ for $i = 1, \ldots, m$.  }
\vskip2mm
\noindent The scaled ABS algorithm for Diophantine equations is now stated.
\vskip2mm

\noindent
{\bf The  ABS Algorithm for  Diophantine Linear Equations}
\begin{enumerate}
\item[(1)] Choose $x_1\in{Z}^n$, arbitrary,  $H_1\in {Z}^{n\times n}$,
      arbitrary  unimodular. Let $i=1$.
\item[(2)] Compute $\tau_i=v_i^Tr_i$ and $s_i=H_iA^Tv_i$.
\item[(3)] {{\it If}} ($s_i=0$ and $\tau_i=0$) {{\it then}} let $x_{i+1}=x_i$,
      $H_{i+1}=H_i$, $r_{i+1}=r_i$ and {{\it go to}} step (5) (the $i$th
      equation is redundant).
      {{\it If}} ($s_i=0$ and $\tau_i\not =0$) {{\it then}} Stop (the
      $i$th equation and hence the system is incompatible).
\item[(4)] $\{s_i\not =0\}$ Compute  $\delta_i=gcd(s_i)$
      and $p_i=H_i^Tz_i$, where
      $z_i\in {Z}^n$ is an arbitrary integer vector satisfying
      $z_i^Ts_i=\delta_i$. {{\it If}} \ \
      $\delta_i\!\!\not| \tau_i$ {{\it then}} Stop 
      (the system is integerly inconsistent),
     {{\it else}}
      Compute
      $
      \alpha_i=\tau_i/\delta_i,
      $
      let
      $
      x_{i+1}=x_i-\alpha_ip_i
      $
 and update $H_i$  by
      $
      H_{i+1}=H_i-\frac{H_iA^Tv_iw_i^TH_i}{w_i^TH_iA^Tv_i}
      $
      where $w_i\in {R}^n$ is an arbitrary integer vector satisfying
       $w_i^Ts_i=\delta_i$.
\item[(5)] {{\it If}} $i=m$ {{\it then}} Stop ($x_{m+1}$ is a solution)
{{\it else}} let $i=i+1$ and {{\it go to}} step (2).
\end{enumerate}

\vspace{0.5cm}

\noindent It follows from Theorem 2 that 
if there exists a solution for  system (15), then
$x=x_{m+1}+H_{m+1}^Tq$, with arbitrary $q\in{Z}^n$, provides all
solutions of (15).  

\vskip8mm
\noindent{\bf  5.  Application of ABS methods to unconstrained optimization}
\vskip5mm

\noindent A direct application of ABS methods to the problem
$$
min f = f(x), \ x \in R^n
$$ 
can be made via the iteration
 $x_{i+1} = x_i - \alpha_i H_i^Tz_i$, where
$H_i$ is reset after $n$ or less steps and $z_i$ is chosen to satisfy 
the descent condition, i.e.
 $g_i^TH_i^Tz_i > 0$, with $g_i$ the gradient of $f(x)$
 in $x_i$. If $f(x)$ is quadratic with Hessian  $A$,
we may identify matrix   $A$ in update (8) with such Hessian.
Defining  
$x' = x_i - \beta v_i$ with $\beta$ a small scalar and noticing that on
quadratic functions one has $g' = g - \beta Av_i$,
update (8) can be written without use of $A$ as

$$H_{i+1} =H_i - H_iy_iw_i^TH_i/w_i^TH_iy_i   \eqno{(16)}$$

with $y_i = g'- g_i$, while the stepsize takes the form
$$
\alpha_i = v_i^Tr_i/p_i^Ty_i \eqno{(17)}.  
$$
The above defined class has termination in at most $n$ steps 
on quadratic functions and on general functions has local convergence
at $n$-step Q-quadratic rate. It is related to a class of so called
projection methods defined by
 Pschenichny and Danilin (1978). Numerical experience with such methods is not
yet available.
\vskip8mm
\noindent {\bf 6. Derivation of general Quasi-Newton updates}
\vskip5mm
 \noindent ABS methods allow to determine the general solution of the
  Quasi-Newton equation, even in presence of additional linear
  conditions as  symmetry and sparsity. Moreover they provide insight
 into the problem of keeping positive definiteness. Notice that while the
 general solution of the
  Quasi-Newton equation was known since a paper of 
 Adachi (1971), the ABS explicit solution for the problem with sparsity and
symmetry was not obtained before.

Let us consider the Quasi-Newton equation defining a new approximation to the
  Jacobian (for a nonlinear system, $g(x)=0$) or to the Hessian
(for an optimization problem, min $f(x)$, with $g(x)$ the gradient
of $f(x)$), here given in the transpose form
$$
d^TB' = y^T
\eqno{(18)}
$$
where $d = x' - x, y = g' - g$. We can notice that (18) is a set of
 $n$ underdetermined linear systems, each one having just one single
 equation and differing only in the right hand side. Thus one step of
 the ABS procedure is enough to get a special solution and a representation
 of all solutions via formula (12). By suitable choice of the parameters
the general solution can be put in the following form, see Spedicato and
Xia (1993)

$$
B' = B - s(B^Td - y)^T/d^Ts + (I - sd^T/d^Ts)Q
\eqno{(19)}
$$
with $Q \in R^{n,n}$  arbitrary and $s \in R^n$ 
arbitrary save that $s^Td \ne 0$. Formula (19) is the same 
 formula derived by  Adachi (1971).

The condition that some elements of  $B'$ are zero or have a constant value
and the condition of symmetry can be written as the following additional 
 linear equations, where
 $b'_i$ is the $i$-th column of $B'$
$$
(b'_i)^Te_j = \eta_{ij}    
$$

\noindent where $\eta_{ij} = 0$ implies sparsity,
$\eta_{ij} = costant$ implies constant value and
  $\eta_{ij} = \eta_{ji}$ implies symmetry.
The above conditions can be dealt with by the ABS algorithm, see
Spedicato and Zhao (1993), obtaining an explicit solution, by successive 
columns in the case of symmetry. By requiring that the 
first $n-1$ diagonal elements
are sufficiently large, then we obtain formulas where
 $B'$ is almost positive definite or almost diagonally dominant (i.e.
 its 
 $(n-1)$-th principal submatrix is positive definite or diagonally dominant).
It is not possible in general to have $B'$ positive definite,
since the last column is completely defined by the
 $n-1$ symmetry conditions and the Quasi-Newton equation.

\vskip8mm

\noindent {\bf 7. The simplex method via the implicit LX algorithm } 
\vskip5mm
\noindent The implicit LX algorithm has a natural application in reformulating the 
simplex method for the LP problem in standard form:
$$ min \ c^Tx$$ subject to
$$Ax = b, \ \ x \ge 0.$$

The suitability of the implicit LX algorithm derives from the fact that
the algorithm, when started with $x_1=0$, generates basic type vectors
  $x_{i+1}$, which are vertices of the polytope defined by the
  constraints of the LP problem if the components not identically zero
  are nonnegative.

The basic step of the simplex method is moving from one vertex to another
according to certain rules and reducing at each step the value of
$c^Tx$.
The direction of  movement  is obtained by solving a linear system,
whose coefficient matrix  $A_B$, the {\it  basic matrix}, is defined by
 $m$ linearly independent columns of   $A$, called the  {\it basic columns}.
This system is usually solved via the LU factorization method or also via
the QR factorization method. 
The new vertex is associated with a new basic matrix
  $A_B'$, obtained by substituting one of the columns of 
 $A_B$ with a column of the matrix
$A_N$ that comprises the columns of $A$ not belonging to
 $A_B$. One has then to solve a new system where just one column has been
 changed. If the LU factorization method is used, then the most efficient
 way to recompute the modified factors is the steepest edge method
 of Forrest and Goldfarb (1992), that requires
 $m^2$ multiplications. This implementation of the simplex method
has the following storage requirement:
$m^2$ for the LU factors and   $mn$ for the matrix
 $A$, that has to be kept to provide the columns needed for the
exchanges.

The reformulation of the simplex method via the implicit LX method, due to 
Xia (1995c), Zhang and Xia (1995), Spedicato, Xia and Zhang (1997), 
 Feng et al. (1997), exploits the fact that in the implicit LX
 method the exchange of the  $j$-th column in $A_B$
with the $k$-th column in $A_N$ corresponds to exchanging previously chosen
parameters
 $z_j = w_j = e_{j_B}$ 
with new parameters  $z_k = w_k =
e_{k_B}$. In term of the implied modification of the Abaffian this operation is a
special case of a general rank-one modification of the parameter matrix
$W=(w_{k_1}, \ldots,w_{k_m})$. The modified Abaffian can be efficiently evaluated
using a general formula of 
 Zhang (1995d), without explicit use of the $k$-th column of
 $A_N$. Moreover all information needed to build the search
 direction (the polytope edge) and to implement the (implicit) column
 exchange is contained in the Abaffian matrix. Thus
 there is no need to keep the matrix $A$. Hence storage requirement is
 only that needed in the construction and keeping of the matrix
  $H_{m+1}$, i.e. respectively
 $n^2/4$ and $n(n-m)$. For   $m$ close to  
 $n$ storage, it is about 8 times less than in the Forrest and Goldfarb implementation
 of the LU method. For small $m$ storage is higher, but there is an
 alternative formulation of the implicit LX method that has a similar
 storage, see Spedicato and Xia (1999).

We now give the main formulas for the simplex method in the classical and
in the ABS formulation. The column in
 $A_N$  that  substitutes a column in
 $A_B$ is usually taken as the column with least relative cost.
In the ABS approach this corresponds to minimize with respect to
  $i \in N_m$ the scalar
$\eta_i = c^TH^Te_i$.
Let $N^\ast$ be the index so chosen. The column in
 $A_B$ to be exchanged is often chosen with the criterion of 
the least edge displacement that keeps nonnegative the basic variables.
Defining
$\omega_i = x^Te_i/e_i^TH^Te_{N^\ast} $
with $x$ the current vertex, then the above criterion is equivalent
to minimize
 $\omega_i$  with respect to the set of indices $i \in B_m$ such that
$$
e_i^TH^Te_{N^\ast} > 0 
\eqno{(20)}
$$
Notice that $H^Te_{N^\ast} \neq 0$ and that an index $i$ such that (20) be
satisfied always exists, unless  $x$ is a solution.

The update of the Abaffian after the exchange of the obtained unit vectors
is given by the following simple rank-one correction
$$
H' = H - (He_{B^\ast} - e_{B^\ast})e_{N^\ast}^{\ T}H/e_{N^\ast}^{\ T}He_{B^\ast}
\eqno{(21)}
$$
The displacement vector $d$, classically obtained by solving system
  $A_Bd = - Ae_{N^\ast}$, is obtained at no cost by
$d = H_{m+1}^{\ T}e_{N^\ast}$. The relative cost vector, classically given by
$r = c - A^TA_B^{-1}c_B$,
with $c_B$ comprising the components of $c$ with indices corresponding to
those of the basic columns, is given by formula
 $r = H_{m+1}c$.

It is easily seen that update (21) requires no more than
 $m(n-m)$ multiplications.  Cost is highest for  $m = n/2$ 
 and gets smaller as  $m$ gets smaller or closer to
  $n$. 
In the  method of  Forrest and
 Goldfarb  $m^2$ multiplications are needed to update the LU
 factors and again $m^2$ to solve the system. Hence the ABS approach via
 formula (21) is faster for $m>n/3$. For $m<n/3$ similar costs are obtained
 using the alternative formulation of the implicit LX algorithm described
 in Spedicato and Xia (1999). For $m$ very close to $n$ the advantage of the
 ABS formulation is essentially of one order and such that 
 no need is seen to develop formulas for the sparse case.

There is an ABS generalization of the simplex method based upon a modification
of the Huang algorithm which is started by a certain singular matrix, see
Zhang (1997). In such generalization the solution is approached via a 
sequence of points lying of the faces of the polytope. If one of such points
happens to be a vertex, then all successive iterates are vertices and
 the standard simplex method is reobtained.

\vskip8mm
\noindent {\bf 8.  ABS unification of feasible direction methods for linearly constrained
minimization}
\vskip5mm
\noindent ABS methods allow a unification of feasible direction methods for linearly
constrained minimization, of which the LP problem is a special case. Let
us first consider the problem with only equality constraints
$$min \ f(x), \ \ x \in R^n$$
subject to
$$ Ax = b, \ \ A \in R^{m,n}, \ \ m \le n,\ \ rank(A) = m.$$

Let $x_1$ be a feasible initial point. If we consider an iteration of the
form
$x_{i+1} = x_i - \alpha_id_i$, then the sequence $x_i$ consists of feasible
points iff
$$
Ad_i = 0  \eqno{(22)}
$$
The general solution of (22) can be written using the ABS formula
(12)
$$
d_i = H_{m+1}^Tq  \eqno{(23)}
$$
In (23) matrix $H_{m+1}$ depends on  parameters 
$H_1$, $W$ and $V$ and $q \in {R^n}$ can also be seen
as a parameter. Hence the general iteration generating feasible points is
$$
x_{i+1} = x_i - \alpha_iH_{m+1}^Tq.  \eqno(24)
$$
The search vector is a descent direction if 
$d^T\nabla{f(x)} = q^TH_{m+1}
\nabla f(x) > 0$. This condition can always be satisfied by a choice of 
 $q$ unless
$H_{m+1}\nabla f(x) = 0$, implying, from the structure of
 $Null(H_{m+1})$, that $\nabla f(x) = A^T\lambda$ for some
$\lambda$, hence that
$x_{i+1}$ is a KT point with $\lambda$ the vector of
 Lagrange multipliers. If
 $x_{i+1}$ is not a KT point, then we can generate descent directions
by taking 
$$q = QH_{m+1}\nabla f(x)   \eqno{(25)}
$$
with $Q$ symmetric positive definite. We obtain therefore a large class
of methods with four parameter matrices ($H_1,W,V,Q$).

Some well-known methods in the literature correspond to taking in (25)
$W=I$ and building the Abaffian as follows.

\begin{description}
\item[(a)]The {\it reduced gradient method} of Wolfe. \ \
 $H_{m+1}$ is built via the implicit LU method.
\item[(b)]The {\it projection method} of Rosen. \ \
 $H_{m+1}$ is built via the Huang
  method.
\item[(c)]The {\it  Goldfarb and Idnani method}. \ \
$H_{m+1}$ is built via a modification of Huang method where
$H_1$ is a symmetric positive definite approximation of the inverse
Hessian of $f(x)$.
\end{description}
To deal with linear inequality constraints there are two approaches in
literature.
\begin{itemize}
\item The {\it active set} method.  
 Here the set of equality constraints is augmented with some inequality constraints,
whose selection varies in the course of the process till the final set
of active constraints is determined. Adding or cancelling a single constraint
corresponds to a rank-one correction to the matrix defining the active
set. The corresponding change in the Abaffian can be performed in order two
operations, see
  Zhang (1995c). For an efficient ABS reformulation of the well-known active
 set method of Goldfarb ande Idnani  see
 Xia, Liu and Zhang (1992) and Xia, Zhang and Liu (1995).
\item The {\it standard form} approach.
Here one uses slack variables to put the problem in the following equivalent
standard form
$$min \ f(x),$$
with the constraints
$$Ax = b, \ \ x \ge 0.$$
\end{itemize}
If $x_1$ satisfies the above constraints, then a sequence of feasible
points is generated by the following iteration
$$
x_{i+1} = x_i - \alpha_i\beta_iH_{m+1}\nabla f(x)    \eqno(26)
$$
where $\alpha_i$ may be chosen by a line search along vector
   $H_{m+1}\nabla f(x)$, while
 $\beta_i > 0$ is chosen to avoid violation of the nonnegativity constraints.

If $f(x)$ is nonlinear, then $H_{m+1}$ can usually be determined once for
all at the first iteration, since generally
 $\nabla f(x)$  changes from point to point, allowing the determination
of a new search direction. However if  $f(x) = c^Tx$ is linear, in which
case we obtain the LP problem, then to get a new search direction we have
to change 
 $H_{m+1}$. We already observed in section 7 that the simplex method
in the ABS formulation is obtained by constructing the matrix
  $H_{m+1}$ via the implicit LX method and at each step modifying one
of the unit vectors used to build the Abaffian. One can show, see
 Xia (1995b), that the  Karmarkar method (to be found in a previous paper by
 Evtushenko (1974)), corresponds to Abaffians built via a variation
 of the Huang algorithm, where the initial matrix is
  $H_1 = Diag(x_i)$ and is changed at every iteration (whether the 
  update of the Abaffian can be performed in order two operations is
  still an open question).  Xia (1995b) discusses several questions
  related to the class of interior feasible point methods obtained via
  such a generalization of Huang algorithm. 
\vskip8mm

\noindent {\bf 9. Other results and conclusion}
\vskip5mm
\noindent Here we recall very briefly a few other results of ABS methods applied
to optimization.
\begin{itemize}
\item The computation of inertia of KT matrices is important in several
algorithms recently proposed for quadratic programming. It can be done
efficiently using a variation of the Huang algorithm, see
Zhang (1999).
\item The LC problem can be reformulated with advantage using the
 implicit LX algorithm, see Feng et al. (1997).
\item ABS methods provide a new computationally interesting approach to  
 concave programming, see Xia and Zhang (1999).
\item   Explicit solution of some classes of integer linear inequalities
and LP problems can be provided via the ABS algorithm for Diophantine equations,
see Esmaeili and Mahdavi-Amiri (2000).
\end{itemize}
Future work on ABS methods will not only deal with theoretical problems where
ABS methods can provide a fresh approach (e.g. the eigenvalue problem or
even combinatorial problems), but will be concentrated on the development
of software, in the framework of the ABSPACK project.

\vskip8mm
\noindent {\bf 10. References}
\vskip5mm
\noindent Here only references are given not quoted in Appendix 2
\vskip5mm
\noindent Adachi N.(1971): On variable metric algorithms, JOTA 7, 391-409\\
Adib M., Mahdavi Amiri N. (2000): ABS-type methods for solving $n$ linear
equations in $n/2$ steps, to appear\\
Dennis J.E. and Turner K. (1987), Generalized conjugate directions, Linear 
Algebra and its Applications 88/89, 187-209\\
Egervary E.(1955), Aufloesung eines homogenen linearen diophantischen
Gleichungsystems mit Hilf von Projektormatrizen, l Publ. Math. Debrecen 4,
481-483\\
Esmaeili M. and Mahdavi-Amiri N. (2000): Solving some linear inequality
systems and LP problems in real and integer spaces via the ABS
algorithm, preprint\\
Evtushenko Y.(1974): Two numerical methods of solving nonlinear
programming problems,  Soviet Doklady Akademii Nauk 251, 420-423\\
Forrest J.J.H. and Goldfarb D. (1992): Steepest edge simplex algorithms for
linear programming, Mathematical Programming  57, 341-174\\
Goldfarb D. and Idnani A. (1983): A numerically stable dual method for
solving strictly convex quadratic programming,
Mathematical Programming 27, 1-33\\
Pschenichny B.N. and Danilin Y.M. (1978):   Numerical Methods in Extremal 
Problems, MIR, Moscow \\
Rosser J.B. (1941): A note on the linear Diophantine equation,  Amer.
Math. Monthly 48, 662-666 \\
Schrijver, A. (1986): Theory of Linear and Integer Programming, John
Wiley and Sons.\\
Spedicato E. and Li Z. (2000): On some classes of iterative ABS-like methods for 
large-scale linear systems, to appear\\
Spedicato E., Bodon E. and Luksan V. (2000), ABSPACK: computational performance of
ABS methods for linear general systems and KT systems, to appear\\
Wedderburn J.H.M. (1934):  Lectures on Matrices, Colloquium Publications,
American Mathematical Society, New York\\
\vskip9mm
{\bf Acknowledgments.} Work partially supported by MURST 1997 program 
Analisi Numerica e Matematica Computazionale.

\newpage
\begin{center}
{\large \bf APPENDIX 1}
\end{center}
\begin{verbatim}   

    RESULTS ON DETERMINED  LINEAR SYSTEMS
 
    Condition  number: 0.21D+20
    IDF2    2000   huang2       0.10D+01  0.69D-11    2000   262.00
    IDF2    2000   mod.huang2   0.14D+01  0.96D-12       4     7.00
    IDF2    2000   lu lapack    0.67D+04  0.18D-11    2000    53.00
    IDF2    2000   qr lapack    0.34D+04  0.92D-12    2000   137.00
    IDF2    2000   gqr lapack   0.10D+01  0.20D-14       3   226.00
    IDF2    2000   lu linpack   0.67D+04  0.18D-11    2000   136.00

    Condition number: 0.10D+61
    IR50    1000   huang2       0.46D+00  0.33D-09    1000    36.00
    IR50    1000   mod.huang2   0.46D+00  0.27D-14     772    61.00
    IR50    1000   lu lapack    0.12D+04  0.12D+04     972     7.00
    IR50    1000   qr lapack    0.63D+02  0.17D-12    1000    17.00
    IR50    1000   gqr lapack   0.46D+00  0.42D-14     772    29.00
    IR50    1000   lu linpack   --- break-down ---                 


    RESULTS ON OVERDETERMINED SYSTEMS
 
    Condition  number: 0.16D+21
    IDF3  1050  950    huang7       0.32D+04  0.52D-13    950    31.00
    IDF3  1050  950    mod.huang7   0.14D+04  0.20D-09      2     0.00
    IDF3  1050  950    qr lapack    0.37D+13  0.83D-02    950    17.00
    IDF3  1050  950    svd lapack   0.10D+01  0.24D-14      2   145.00
    IDF3  1050  950    gqr lapack   0.10D+01  0.22D-14      2    27.00

    Condition number:  0.63D+19
    IDF3  2000  400    huang7       0.38D+04  0.35D-12    400     9.00
    IDF3  2000  400    mod.huang7   0.44D+03  0.67D-12      2     0.00
    IDF3  2000  400    impl.qr5     0.44D+03  0.62D-16      2     0.00
    IDF3  2000  400    expl.qr      0.10D+01  0.62D-03      2     0.00
    IDF3  2000  400    qr lapack    0.45D+12  0.24D-02    400     8.00
    IDF3  2000  400    svd lapack   0.10D+01  0.65D-15      2    17.00
    IDF3  2000  400    gqr lapack   0.10D+01  0.19D-14      2    12.00
    

    RESULTS ON UNDERDETERMINED LINEAR SYSTEMS

    Condition number:  0.29D+18
    IDF2   400 2000    huang2       0.12D-10  0.10D-12    400    12.00
    IDF2   400 2000    mod.huang2   0.36D-08  0.61D-10      3     1.00
    IDF2   400 2000    qr lapack    0.29D+03  0.37D-14    400     9.00
    IDF2   400 2000    svd lapack   0.43D-13  0.22D-14      3    68.00
    IDF2   400 2000    gqr lapack   0.18D-13  0.24D-14      3    12.00
  
    Condition number:  0.24D+19
    IDF3   950 1050    huang2       0.00D+00  0.00D+00    950    33.00
    IDF3   950 1050    mod.huang2   0.00D+00  0.00D+00      2     1.00
    IDF3   950 1050    qr lapack    0.24D+03  0.56D-14    950    17.00
    IDF3   950 1050    svd lapack   0.17D-14  0.92D-16      2   178.00
    IDF3   950 1050    gqr lapack   0.21D-14  0.55D-15      2    26.00
	
    RESULTS ON KT SYSTEMS

    Condition number:  0.26D+21 
    IDF2  1000  900    mod.huang    0.55D+01  0.23D-14     16    24.00
    IDF2  1000  900    impl.lu8     0.44D+13  0.21D-03   1900    18.00
    IDF2  1000  900    impl.lu9     0.12D+15  0.80D-02   1900    21.00
    IDF2  1000  900    lu lapack    0.25D+03  0.31D-13   1900    62.00
    IDF2  1000  900    range space  0.16D+05  0.14D-11   1900    87.00
    IDF2  1000  900    null space   0.89D+03  0.15D-12   1900    93.00
 
    Condition number:  0.70D+20 
    IDF2  1200  600    mod.huang    0.62D+01  0.20D-14     17    36.00
    IDF2  1200  600    impl.lu8     0.22D+07  0.10D-08   1800    44.00
    IDF2  1200  600    impl.lu9     0.21D+06  0.56D-09   1800    33.00
    IDF2  1200  600    lu lapack    0.10D+03  0.79D-14   1800    47.00
    IDF2  1200  600    range space  0.11D+05  0.15D-11   1800    63.00
    IDF2  1200  600    null space   0.38D+04  0.13D-12   1800   105.00
	\end{verbatim}

\newpage
\begin{center}
{\large \bf APPENDIX 2}
\end{center}
\vskip5mm

\noindent Vespucci et al. (1992) have given a 205 entries bibliography of ABS methods 
up to the end of 1991, Nicolai and Spedicato (1997) have updated that
bibliography to just over  300 entries. The present bibliography,  
tentatively complete up to the end of 1999, contains 349 entries (not
counting earlier versions of some papers).

The following abbreviations are used:
\vskip2mm
\noindent
JOTA \mbox{}\hspace{2cm} Journal of Optimization Theory and Applications\\
OMS \mbox{}\hspace{2.2cm} Optimization Methods and Software \\
ABSC1 \mbox{}\hspace{1.8cm} Proceedings of the First International Conference
on ABS\\
\mbox{}\hspace{3.2cm} Algorithms, Luoyang, September 1991,\\
\mbox{}\hspace{3.2cm} University of Bergamo (Spedicato E. editor), 1992\\
ABSC2 \mbox{}\hspace{1.8cm} Proceedings of  the Second International\\
\mbox{}\hspace{3.3cm}Conference on ABS Algorithms,                    \\
\mbox{}\hspace{3.3cm}Beijing, June 1995, Report  DMSIA 19/97,           \\
\mbox{}\hspace{3.3cm}University of Bergamo                                \\
\mbox{}\hspace{3.3cm}(Spedicato E., Li Z. and Luksan L. editors)\\
ABSDAL \mbox{}\hspace{1.4cm} Collection of reports on the ABS algorithms,\\
\mbox{}\hspace{3.2cm} Dalian University of Technology, 1991\\
QDMSIA \mbox{}\hspace{1.4cm} Quaderni del Dipartimento di Matematica,
Statistica,\\
\mbox{}\hspace{3.2cm} Informatica  ed  Applicazioni, University of Bergamo\\
\vskip5mm
Abaffy J. (1979), {\sl A \ \  line\'{a}ris \ \  egyenletrendszerek \ \ 
\'{a}ltal\'{a}nos \ \  megold\'{a}s\'{a}nak \ \  egy \ \  
m\'{o}dszeroszt\'{a}lya}, Alkalmazott Matematikai Lapok 5, 223-240.

Abaffy J.(1984), {\sl A generalization of the ABS class}, Proceedings
of the Second Conference  on Automata Languages  and Mathematical Systems, 
7-11, University of Economics, Budapest.

Abaffy J. (1986), {\sl Some special cases of the  ABS  class for band
matrices},  Proceedings of the Fourth Conference on Automata Languages 
and Mathematical Systems, 9-15, University of Economics, Budapest.

Abaffy J. (1987a), {\sl Error analysis and stability for some cases in 
the ABSg class}, Report 190, Numerical Optimisation Centre, 
Hatfield Polytechnic.

Abaffy J. (1987b), {\sl The description of a package for sparse linear 
systems  using six algorithms  in the ABSg class},  Report 192,  
Numerical Optimisation Centre, Hatfield Polytechnic.

Abaffy J. (1987c),  {\sl Preliminary test results with some algorithms of 
the ABSg Class},  Report 193, Numerical Optimisation Centre, 
Hatfield Polytechnic.

Abaffy J. (1987d),  {\sl Optimal methods in Broyden's stability sense}, 
Report 194, Numerical Optimisation Centre, Hatfield Polytechnic.

Abaffy J. (1988a),  {\sl A superlinear convergency theorem in the ABSg class 
for nonlinear algebraic equations}, JOTA 59, 39-43.

Abaffy J. (1988b),  {\sl Equivalence of a generalization  of Sloboda's 
algorithm  with a subclass  of the generalized ABS algorithm  for linear 
systems}, QDMSIA 1/88.

Abaffy J. (1988c),  {\sl Derivation of Lanczos type methods in the 
generalized  ABS algorithm for linear systems},  QDMSIA  2/88.

Abaffy J. (1988d),  {\sl Formulation of the  conjugate gradient  type 
algorithm of Fridman for indefinite symmetric systems in the subclass of 
the generalized symmetric ABS algorithms}, QDMSIA 3/88.

Abaffy J. (1988e),  {\sl Reorthogonalized  methods  in the ABSg class},
QDMSIA 17/88.

Abaffy J. (1991),  {\sl ABS algorithms for sparse linear systems}, 
in Computer Algorithms  for Solving Linear Algebraic Equations: 
the State of the Art, NATO ASI Series,  Vol. F77,  111-132, 
Springer-Verlag (Spedicato E. editor).

Abaffy J. (1995), {\sl A new  method for computing eigenvalues}, ABSC2.
Abaffy J. and Bertocchi M. (1990), {\sl Vectorized ABS conjugate direction 
algorithms  on the IBM 3090 VF},  Pure Mathematics and Applications, Ser.B, 
1, 189-202.

Abaffy J. and Bertocchi M. (1991),  {\sl A vectorized FORTRAN code  of the 
implicit conjugate direction algorithm of the ABS class}, Manual 1/57, 
IAC, Rome.

Abaffy J., Bertocchi M. and Torriero A. (1991a), {\sl Characterization of 
M-matrices via ABS methods  and its application to input-output analysis}, 
QDMSIA 11/91.

Abaffy J., Bertocchi M. and Torriero A. (1991b),  {\sl Una caratterizzazione 
delle M-matrici attraverso la classe ABS: una applicazione al modello 
input-output}, Atti del XV Convegno dell'Associazione  per la Matematica 
Applicata  alle Scienze Economiche e Sociali, Grado, September 1991. 

Abaffy J., Bertocchi M. and Torriero A. (1992a), {\sl Criteria for
transforming a Z-matrix into an M-matrix}, OMS 1, 183-196. Also as QDMSIA 12/91.

Abaffy J., Bertocchi M. and Torriero A. (1992b), {\sl Analysis of a subclass 
of nonsingular M-matrices and its economic applications}, QDMSIA 22/92.

Abaffy J., Bertocchi M. and Torriero A. (1992c), {\sl Un algoritmo per la 
determinazione dell' autovalore di Frobenius di matrici non negative}, Atti 
XVI Convegno AMASES, Treviso, 763-767.

Abaffy J., Bertocchi M. and Torriero A. (1992d), {\sl Perturbations of M-matrices 
via ABS methods and their applications to input-output analysis}, QDMSIA 23/92.
Abaffy J., Bertocchi M. and Torriero A. (1992e), {\sl New algorithms for computing 
the Frobenius eigenvalue}, QDMSIA 25/92.

Abaffy J., Bertocchi M. and Torriero A. (1994), {\sl Sul calcolo dell' 
autovalore di Frobenius di matrici non negative}, Scritti in onore di Giovanni 
Melzi, Vita e Pensiero, Milano, 9-26 (Manara C.F. and Faliva M. editors).

Abaffy J.,  Broyden C.G.  and  Spedicato E.  (1982), {\sl A  class of direct 
methods  for linear systems II:  generalization, non-singular  representation 
and  other tales}, Report SOFMAT 21/82,  IAC, Rome.

Abaffy J., Broyden C.G. and Spedicato E. (1983),  {\sl Numerical  performance  
of the pseudosymmetric algorithm  in the ABS class versus LU  factorization 
with iterative refinement}, Report SOFMAT 8/83, IAC, Rome.

Abaffy J., Broyden C.G. and Spedicato E. (1984a), {\sl A class  of direct 
methods of quasi-Newton type for general linear systems}, Monografia di Software 
Matematico 33, IAC, Rome. 

Abaffy J., Broyden C.G. and Spedicato E. (1984b), {\sl A  class   of direct 
methods  for linear systems}, Numerische Mathematik, 45, 361-376.

Abaffy J., Broyden C.G. and Spedicato E. (1986), {\sl Conjugate direction 
methods  for linear and nonlinear systems  of algebraic equations}, 
QDMSIA 17/86.

Abaffy J. and Dixon L.C.W. (1987), {\sl On solving sparse band systems with 
three algorithms  of the ABS family}, Report 191, Numerical 
Optimisation Centre, Hatfield Polytechnic.

Abaffy J. and Gal\'{a}ntai A. (1987), {\sl Conjugate  direction  methods for 
linear and nonlinear systems of algebraic equations}, Colloquia Mathematica 
Societatis J\'{a}nos Bolyai 50,  Numerical Methods, Miskolc, 481-502, 
North-Holland (Stoyan G. editor).

Abaffy J.,  Gal\'{a}ntai A.  and  Spedicato E.  (1984), {\sl Convergence 
properties of the ABS algorithm for nonlinear algebraic equations}, 
QDMSIA 7/84.

Abaffy J.,  Gal\'{a}ntai A.  and  Spedicato E.  (1987a), {\sl The  local 
convergence of ABS methods for non linear algebraic equations}, Numerische 
Mathematik, 51, 429-439. Also as QDMSIA 16/85

Abaffy J., Gal\'{a}ntai A. and Spedicato E. (1987b), {\sl Application  of the 
ABS class  to  unconstrained function minimization}, QDMSIA 14/87.

Abaffy J.,  Gal\'{a}ntai A.  and  Spedicato E. (1989), {\sl Forward error 
analysis in ABS methods}, Preprint, University of Bergamo.

Abaffy J. and Spedicato E. (1983a), {\sl A  code for general linear systems},  
Monografia di Software Matematico 21, IAC, Rome.

Abaffy J. and Spedicato E. (1983b), {\sl Computational experience with a 
class of direct methods for linear systems},  Report SOFMAT 1/83, IAC, Rome.

Abaffy J. and Spedicato E. (1983c), {\sl On the symmetric algorithm in the ABS 
class of direct methods  for linear systems}, Report SOFMAT 7/83, 
IAC, Rome.

Abaffy J. and Spedicato E. (1984), {\sl A generalization of Huang's method  
for solving systems of linear algebraic equations}, Bollettino Unione 
Matematica Italiana 3-B, 517-529.

Abaffy J. and Spedicato E. (1985), {\sl A generalization of the ABS algorithm for linear
systems}, Computers and Computing, Proceedings of the
International Conference on Future Trends of Computing, Grenoble, December 
1985, Wiley and Sons (Chenin P., Di Crescenzo C. and Robert F. editors).
Also as QDMSIA 4/85.

Abaffy J. and Spedicato E. (1987), {\sl Numerical experiments with the 
symmetric algorithm in the ABS class for linear systems}, Optimization, 
18, 197-212. Also as QDMSIA 18/85.  

Abaffy J. and Spedicato E. (1988), {\sl Orthogonally scaled and optimally 
stable algorithms in the scaled ABS class  for linear systems}, QDMSIA 8/88.

Abaffy J. and Spedicato E. (1989a), {\sl ABS  Projection Algorithms: 
Mathematical Techniques for  Linear and Nonlinear Algebraic Equations}, 
Ellis Horwood, Chichester (Chinese translation published in 1991 by Beijing
Polytechnical University Press, Beijing. Russian
translation with update appendix published in 1996 by Mir, Moscow).

Abaffy J. and Spedicato E. (1989b), {\sl On the use  of the ABS algorithm for 
some linear programming problems}, Preprint, University of Bergamo.

Abaffy J. and Spedicato E. (1990), {\sl A  class of scaled direct methods for 
linear systems}, Annals of the Institute of Statistical Mathematics, 42, 
187-201.

Adib M, Mahdavi Amiri N. and Spedicato E. (1999), {\sl Broyden's method as
and ABS algorithm}, Publ. Univ. Miskolc, Series D, 40, 3-13.

Belyankov A.Y., (1990), {\sl Exhaustion methods in numerical linear algebra},
Proceedings on International Conference on Optimization Software, Birkhauser 

Bertocchi M. (1989), {\sl Numerical experiments  with ABS algorithms for 
linear systems on a parallel machine}, JOTA 60, 375-392.

Bertocchi M. (1991), {\sl A generalization of the orthogonally and optimally 
stable scaled ABS algorithms for linear systems}, Pure Mathematics and 
Applications, Serie B, 2, 251-266.

Bertocchi M. and Spedicato E. (1988), {\sl Vectorizing the modified Huang algorithm of the ABS class on the IBM 3090 VF}, QDMSIA 23/88. Also in 
Proceedings of the International Meeting on Parallel Computing, Verona, 
Italy, 1988.

Bertocchi M. and Spedicato E. (1989a), {\sl \ Performance \ of \ the \ implicit 
Gauss- Cholesky algorithm of the ABS class on the IBM 3090 VF}, Proceedings of 
the 10-th Symposium on Algorithms, Strbske Pleso, 1989, 30-34.

Bertocchi M. and Spedicato E. (1989b), {\sl Computational performance on the
IBM 3090 VF of the  modified  Huang and the implicit Gauss-Cholesky algorithms  
versus the Gaussian solver in the ESSL library  on ill-conditioned problems},   
Proceedings of the 
10-th Symposium on Algorithms,  Strbske Pleso, April 17-21, 1989, 22-29.
Also as QDMSIA 14/88.

Bertocchi M. and Spedicato E. (1989c), {\sl Vectorizing the implicit 
Gauss-Cholesky algorithm of the ABS class on the IBM 3090 VF}, 
Report 1/37, IAC, Rome. Also as QDMSIA 22/88. 

Bertocchi M. and Spedicato E. (1990a), {\sl Block ABS algorithms for dense  
linear  systems in a  vector  processor  environment},  
Proceedings of the Conference on Supercomputing Tools for Science and 
Engineering, Pisa, December 1989, Franco Angeli, Milano, 39-46 
(Laforenza D. and Perego R. editors). Also as QDMSIA 2/90.

Bertocchi M. and Spedicato E. (1990b), {\sl A vectorized FORTRAN code of the 
implicit Gauss-Cholesky  algorithm  of the ABS class}, Manual 1/34,  IAC, Rome.

Bertocchi M., Brandolini L. and Spedicato E. (1990a), {\sl A vectorized 
FORTRAN code of the modified Huang algorithm  of the ABS class}, Manual 1/35, 
IAC, Rome.

Bertocchi M., Brandolini L. and Spedicato E. (1990b), {\sl Vectorizing the 
modified Huang algorithm  of the ABS class  on the IBM 3090VF}, 
Report 1/36, IAC, Rome.

Bertocchi M., Spedicato E. and Vespucci M.T. (1991), {\sl A vectorized FORTRAN 
code of the implicit QR algorithm  of the ABS class}, Manual 1/56,  IAC, Rome.

Bertocchi M., Spedicato E. and Vespucci M.T. (1992), {\sl Vectorizing the 
implicit QR algorithm  of the ABS class  on the IBM 3090 VF}, 
Proceedings of the Meeting on Parallel 
Computing: Problems, Methods and Applications, Elsevier, 
99-108 (Messina P. and  Murli A. editors).

Bodon E. (1987), {\sl Algoritmi ABS per minimi quadrati ed equazioni non 
lineari}, Ph.D. Dissertation, University of Genova.

Bodon E. (1989a), {\sl A  code  for linear full rank least squares problems 
based upon an implicit QR ABS algorithm}, QDMSIA 14/89.

Bodon E. (1989b), {\sl Globally  optimally  conditioned update in the ABS 
class for linear systems}, QDMSIA 17/89.

Bodon E. (1989c), {\sl Biconjugate  algorithms  in the  ABS  class I: 
alternative formulation and theoretical properties}, QDMSIA 3/89.

Bodon E. (1990a), {\sl A code for solving determined or underdetermined 
full or deficient rank linear systems based upon the optimally conditioned 
ABS algorithm},  Preprint,  University of Bergamo.

Bodon E. (1990b), {\sl Numerical behaviour of the implicit QR algorithm
of the ABS class for nonlinear systems}, Preprint, University of Bergamo.

Bodon E. (1992a), {\sl A code of the Hestenes-Stiefel type ABS algorithm}, 
QDMSIA 16/92.

Bodon E. (1992b), {\sl Numerical experiments with Gauss-ABS algorithms
on tridiagonal systems of linear equations}, QDMSIA 31/92.

Bodon E. (1993a), {\sl Numerical results on the ABS algorithms  for
linear systems of equations}, QDMSIA 9/93.

Bodon E. (1993b), {\sl Numerical  experiments  with ABS algorithms on upper banded 
systems of linear equations}, Report TR/PA93/14, 
CERFACS, Toulouse. Also as QDMSIA 17/92.

Bodon E. (1993c), {\sl Numerical  experiments  with ABS algorithms on 
banded systems of linear equations}, Report TR/PA/93/13, 
CERFACS, Toulouse. Also as QDMSIA 18/92.

Bodon E. (1993d), {\sl Computational performance of the implicit
QR algorithm for linear least squares}, QDMSIA 17/93.

Bodon E. (1993e), {\sl Two codes of block ABS algorithms for linear
systems of equations}, Preprint, University of Bergamo.

Bodon E. (1995a), {\sl Some experiments with the ABS-GMRES method},
Preprint, University of Bergamo.

Bodon E. (1995b), {\sl Two codes of the ABS-GMRES method}, Preprint,
University of Bergamo.

Bodon E. and Spedicato E. (1990a), {\sl The scaled  ABS algorithm and its 
application  to  linear  least  squares:  some numerical results}, 
QDMSIA 3/90.

Bodon E. and Spedicato E. (1990b), {\sl Numerical evaluation of  the implicit 
LU, LQ and QU algorithms in the ABS class}, QDMSIA 28/90.

Bodon E. and Spedicato E. (1991), {\sl Factorized ABS algorithms for linear 
systems:  derivation and numerical results}, QDMSIA 14/91.

Bodon E. and Spedicato E. (1992), {\sl On some STOD-ABS algorithms for large 
linear systems},  QDMSIA 15/92.

Broyden C.G. (1985), {\sl On the numerical stability of Huang's and related
methods}, JOTA 47, 401-412.

Broyden C.G. (1991), {\sl On the numerical stability of Huang's update}, 
Calcolo, 28 (1991), 303-311. Also as QDMSIA 18/89.

Ceribelli C. (1989), {\sl Implementazione dell' algoritmo di Huang per sistemi 
lineari a bande}, Ms. Dissertation, University of Bergamo.

Chen F. and Xia Z. (1999), {\sl Roundoff error analysis of the ABS algorithm},
QDMSIA 10/99

Chen Z. and  Deng  N. (1992), {\sl A general algorithm for underdetermined 
linear systems}, ABSC1.

Chen Z. and Deng N. (1993), {\sl An algorithm with $n$-step quadratic convergency
for nonlinear systems}, Numerical Mathematics, a  Journal of Chinese Universities
15, 32-39

Chen Z. and Deng N. (1995), {\sl A class of algorithms for quadratic
programming via the ABS approach}, ABSC2.

Chen Z. and Deng N. (1997),{\sl Some algorithms for the convex quadratic
programming problem via the ABS approach}, OMS 8, 157-170

Chen Z., Deng N. and Spedicato E. (1999), {\sl Truncated difference ABS-type algorithm
for nonlinear systems of equations}, OMS 10, 1-12

Chen Z. and Spedicato E. (1993), {\sl Computing the null space of a matrix 
via the ABS approach}, Preprint, University of Bergamo.

Chen Z. and Xue Y. (1992), {\sl Efficient ABS-type algorithms for solving
linear systems of equations}, Preprint, Polytechnical University, Beijing.

Deng N. (1998). {\sl A new nonlinear ABS method with convergence and
efficiency analysis}, QDMSIA 3/98.

Deng N. and Chen Z. (1990), {\sl Truncated nonlinear ABS algorithm
and its convergence property}, Computing 45, 169-173.

Deng N. and Chen Z. (1991), {\sl A  unified  convergence theorem for nonlinear 
ABS-type algorithm},  Preprint, University of Bergamo.

Deng N. and Chen Z. (1998), {\sl A new nonlinear ABS-type algorithm and its
efficiency analysis}, OMS 10, 71-86

Deng N. and Spedicato E. (1988), {\sl Optimal conditioning parameter selection 
in the ABS class  through a rank two update formulation}, QDMSIA 18/88.

Deng N. and Xiao Y. (1989), {\sl The minimum norm correction class for linear 
system of equations}, Preprint, Department of Applied Mathematics, Beijing 
Polytechnic University.

Deng N. and Wang L. (1993), {\sl A modified ABS-type algorithm depending on 
nonlinearity}, Journal of Beijing Agricultural Engineering University 13, 26-32.

Deng N. and Zhu M. (1992a),  {\sl Further results on the  local
convergence of the nonlinear ABS algorithm}, OMS 1, 211-222. Also as QDMSIA 15/91.

Deng N. and Zhu M. (1992b), {\sl The local convergence of the nonlinear 
Voevodin-ABS algorithm}, OMS 1, 223-231. Also as QDMSIA 7/91.

Deng N., Li Z. and Spedicato E. (1996), {\sl On sparse Quasi-Newton 
quasi-diagonally dominant updates}, QDMSIA 1/96.

Deng N., Spedicato E. and Vespucci M.T. (1991), {\sl Experiments with the ABS 
implicit Gauss-Cholesky algorithm  on nested dissection matrices},  
Report 1/69, IAC, Rome.

Deng N., Spedicato E. and Zhu M. (1994), {\sl A local convergence
theorem for nonlinear ABS algorithms}, Computational  and
Applied Mathematics 13, 49-59. Also as QDMSIA 9/91.

Deng N., Xiao Y. and Zhou F. (1992), {\sl Variationally  derived  ABS methods}, 
OMS 1, 197-210.

Dixon L.C.W. (1996), {\sl On the numerical stability of the implicit LU
method and related methods}, QDMSIA 18/97.

Egervary E. (1953), {\sl On a property of the projector matrices and its
application to canonical representation of matrix functions}, Acta
Scientiarum Mathematicarum 15, 1-6, 1953.

Egervary E. (1957), {\sl \"{U}ber eine Methode zur numerische L\"{o}sung der
Poissonischen Differentialgleichung f\"{u}r beliebige Gebiete}, Acta Mathematica
Academiae Scientiae Hungaricae 9, 341-361.

Egervary E. (1960), {\sl On rank-diminishing operations and their
application to the solution of linear equations}, ZAMP 9, 377-386.

Esmaeili H. (1998), {\sl Application of the ABS algorithm to Diophantine
equations and other tales}, QDMSIA 4/98.

Esmaeili H., N. Amiri-Mahdavi and Spedicato E. (1999), {\sl Solution of
Diophantine systems via the ABS methods}, QDMSIA 29/99

Eremin Y. (1991), {\sl A block generalization of the ABS algorithm}, Preprint,
Institute of Numerical Analysis, Russian Academy of Science, Moscow.

Faridpour H. (1994), {\sl A class of scaled direct methods for solving a 
system of linear equations}, Ms. Dissertation, Tabriz University. 

Fegyverneki S. (1995), {\sl Some new estimators for robust nonlinear 
regression}, Publications of University of Miskolc, Series D, 
36, 23-32.

Feng E. (1994), {\sl Implementing and testing an ABS algorithm for
linearly constrained optimization}, QDSMIA 17/94.

Feng E. (1995), {\sl The FORTRAN code  LCABSM of an ABS algorithm
for linearly constrained optimization}, QDSMIA 1/95.

Feng E. and Wang X. (1991), {\sl Application of the ABS algorithm to solve 
the mathematical simulation  of geological history}, ABSDAL.

Feng E., Shi G. and Zhang L. (1994), {\sl The ABS formulation of weighted
generalized inverse}, Northeast Operations Research, 8, 39-41.

Feng E., Wang X. and Wang X.L. (1995), {\sl The application of ABS algorithm
in the dual simplex method}, ABSC2.

Feng E., Wang X. and Wang X.L. (1997), {\sl On the application of the ABS
algorithm to linear programming and linear complementarity}, OMS 8, 133-142

Fodor S. (1999a), A class of ABS methods for Diophantine systems of equations,
part I, preprint, University of Economics, Budapest

Fodor S. (1999b), A class of ABS methods for Diophantine systems of equations,
part II, preprint, University of Economics, Budapest

Fragnelli V. and Resta G. (1985), {\sl An  ABS method for solving suitably 
structured linear systems  is supported by parallel architectures}, Proceedings  
of the First European Workshop on Parallel Processing Techniques for 
Simulation, Manchester, Plenum Press, 234-247. 

Gal\'{a}ntai A. (1989), {\sl A new convergence theorem for nonlinear ABS 
methods}, QDMSIA 12/89.

Gal\'{a}ntai A. (1990a), {\sl Block ABS methods  for nonlinear systems of 
algebraic equations}, in Computational Solution of Nonlinear Systems of Equations, 
SIAM Lectures in Applied Mathematics 26, 181-187 
(Allgower E.L. and Georg K. editors).

Gal\'{a}ntai A. (1990b),  {\sl A global convergence theorem for nonlinear ABS 
methods}, QDMSIA 9/90.

Gal\'{a}ntai A. (1991a), {\sl Convergence theorems  for the nonlinear ABS 
methods}, Colloquia Mathematica Societatis J\'{a}nos Bolyai 59, Numerical 
Methods, Miskolc 1990, North-Holland, 65-79 (Stoyan G. editor). 

Gal\'{a}ntai A. (1991b), {\sl Analysis  of error propagation in the ABS class}, 
Annals of the Institute of Statistical Mathematics 43, 597-603.
Also as QDMSIA 1//87.

Gal\'{a}ntai A. (1993a), {\sl Generalized implicit LU algorithms in the
class of  ABS methods for linear and nonlinear systems of
algebraic equations}, QDMSIA 5/93.

Gal\'{a}ntai A. (1993b), {\sl ABS methods on large nonlinear systems with banded 
Jacobian structure}, QDMSIA 14/93.

Gal\'{a}ntai A. (1993c), {\sl Testing of implicit LU ABS methods on large 
nonlinear systems with banded Jacobian}, QDMSIA 19/93.

Gal\'{a}ntai A. (1994), {\sl ABS methods for the solution of linear and 
nonlinear systems of equations}, Habilitation Dissertation, 
Technical University of Budapest.

Gal\'{a}ntai A. (1995a), {\sl Hybrid Quasi-Newton methods for nonlinear
equations}, ABSC2.

Gal\'{a}ntai A. (1995b), {\sl The global convergence of the ABS methods for
a class of nonlinear problems}, OMS 4, 283-295. Also as QDMSIA 9/92.

Gal\'{a}ntai A. (1995c), {\sl A fast method for solving systems of nonlinear 
algebraic equations}, Publications of the University of Miskolc, Series C, 
45, 85-92.

Gal\'{a}ntai A. (1997a), {\sl Hybrid Quasi-Newton ABS methods for nonlinear
equations}, QDMSIA 11/97.

Gal\'{a}ntai A. (1997b), {\sl The local convergence of Quasi-Newton ABS
methods},  Publications of the University of Miskolc, Series D, 37, 31-38

Gal\'{a}ntai A. (1999a), {\sl Perturbation theory for full rank factorizations},
QDMSIA 40/99

Gal\'{a}ntai A. (1999b), {\sl Parallel ABS projection methods for linear and
nonlinear systems with block arrohead structure}, Computers and Mathematics
with Applications 38, 11-17

Gal\'{a}ntai A. and Jeney A. (1993), {\sl Quasi-Newton ABS  methods},
Proceedings of MICROCAD-93 Conference, Miskolc, 63-68.

Gal\'{a}ntai A. and Jeney A. (1996), {\sl Quasi-Newton ABS methods for solving 
nonlinear algebraic systems of equations}, JOTA 89, 561-573.

Gal\'{a}ntai A., Jeney A. and Spedicato E. (1993a), {\sl Testing  of
ABS-Huang methods on moderately large nonlinear systems}, QDMSIA 6/93

Gal\'{a}ntai A., Jeney A. and Spedicato E. (1993b), {\sl A FORTRAN code
of a Quasi-Newton ABS algorithm}, QDMSIA 7/93.

Ge R. (1997), {\sl A family of Broyden-ABS type methods for solving nonlinear
systems}, QDMSIA 1/97.

Ge R. and Xia Z. (1996), {\sl An ABS algorithm for solving nonsingular
systems}, QDMSIA 36/96.

Gecan A. and Snider A.D. (1997), {\sl Rationale for the Camerini-Fratta-Maffioli
modified subgradient deflection and its ABS formulation}, ABSC2, 57-66

Gredzhuk V.I. and Petrina D.Y. (1998), {\sl The ABS method and the method
of averaging functional corrections}, QDMSIA 8/98

Giudici S. (1993), {\sl Applicazione del metodo ABS per la soluzione
dei problemi di minimi quadrati nonlineari}, Ms. Dissertation, University of Bergamo.

Gu S. (1989), {\sl The simplified ABS algorithms}, QDMSIA 5/89.

Huang H.Y. (1975), {\sl A direct method for the general solution of a system 
of linear equations}, JOTA 16, 429-445.

Huang Z. (1989), {\sl On the  ABS methods for nonlinear problems}, 
Ph.D. Dissertation, Fudan University, Shanghai.

Huang Z. (1991a), {\sl Multi-step  nonlinear  ABS methods and their efficiency 
analysis}, Computing 46, 143-153. 

Huang Z. (1991b), {\sl Modified  ABS methods for nonlinear systems without 
evaluating Jacobian matrices}, Numerical Mathematics, a Journal of Chinese Universities 
13, 60-71.

Huang Z. (1992a), {\sl Convergence analysis of the nonlinear block scaled
ABS methods}, JOTA 75, 331-344. Also as QDMSIA 14/90.

Huang Z. (1992b), {\sl A family of discrete ABS type algorithms for solving 
systems of nonlinear equations}, Numerical Mathematics, a Journal of Chinese 
Universities 14, 130-143.

Huang Z. (1992c), {\sl ABS methods for solving nonlinear least squares problems}, 
Communications on Applied Mathematics and Computation 6, 75-85.

Huang Z. (1993a), {\sl A class of ABS methods for nonlinear least squares
problems}, Ricerca Operativa 26, 67-78. Also as QDMSIA 3/92.

Huang Z. (1993b), {\sl On the convergence analysis of the nonlinear ABS methods}, 
Chinese Annals of Mathematics 14 B, 213-224.

Huang Z. (1993c), {\sl Row  update ABS  methods for 
solving sparse nonlinear systems of equations}, 
OMS 2, 297-309.

Huang Z. (1993d), {\sl Restart row update ABS methods for solving systems of 
nonlinear equations}, Computing 50, 229-239. Also as QDMSIA 13/92.

Huang Z. (1994a), {\sl Row update ABS methods for systems of nonlinear
equations}, OMS 3, 41-60. Also as QDMSIA 14/92.

Huang Z. (1994b), {\sl Generalization of row update ABS methods
for systems of nonlinear equations}, QDMSIA 14/94.

Huang Z. (1994c), {\sl A new method for solving nonlinear underdetermined
systems}, Computational and Applied Mathematics 1, 33-48. Also as QDMSIA 5/92.

Huang Z. and Spedicato E. (1994d), {\sl Numerical experiments with
Newton's, Brent's and Huang type methods for nonlinear
systems}, QDMSIA 1/94.

Huang Z. and Spedicato E. (1994e), {\sl Numerical testing of Quasi-Newton
and some related methods for nonlinear systems}, QDMSIA 4/94.

Huang Z. and Spedicato E. (1994f), {\sl Further numerical testing
of Newton-type methods for nonlinear systems}, QDMSIA 9/94.

Huang Z. and Spedicato E. (1996), {\sl ABS formulation of the GMRES method},
QDMSIA 14/96.

Jeney A. (1991), {\sl Numerical studies of ABS methods for solving 
nonlinear systems of equations}, Alkalmazott Matematikai Lapok, 
15, 353-364 (in Hungarian).

Jeney A. (1991), {\sl Discretization of a subclass of ABS methods for 
solving nonlinear systems of equations}, Alkalmazott Matematikai 
Lapok, 15, 365-379 (in Hungarian). Also as QDMSIA 10/92 (in English).

Jeney A. (1992a), {\sl The performance of ABS-methods 
for solving nonlinear systems of algebraic equations}, QDMSIA 11/92.

Jeney A. (1992b), {\sl Discrete ABS methods for nonlinear systems of equations}, 
Proceedings of MICROCAD-92, Miskolc, 107-110 (in Hungarian).

Jeney A. (1995), {\sl The convergence of discrete ABS methods and numerical
investigations}, Proceedings of MICROCAD-95, Miskolc, Section K, 22-25.

Jeney A. (1996a), {\sl Calculation of deformations of fibre-reinforced hoses 
with system of nonlinear algebraic equations}, Proceedings of MICROCAD-96, 
Section L, Miskolc, 25-30.


Jeney A. (1996c), {\sl Numerical comparison of multistep Quasi-Newton
ABS methods with the multistep Broyden method for solving systems of
nonlinear equations}, Publications of the University of Miskolc, Series D,
36, 47-54

Jeney A. (1997a), {\sl The local convergence of Quasi-Newton ABS methods
}, Publications of the University of Miskolc, 
Series D, 37, 31-38.

Jeney A. (1997b), {\sl Multistep discrete ABS methods for solving systems
of nonlinear equations}, Publications of the University of Miskolc, 
Series D, 37, 39-46.

Lai Y., Gao Z. and He G. (1995), {\sl An ABS-generalized projection algorithm
for linearly constrained optimization problems}, ABSC2.

Li C. (1997), {\sl On the preconditioning properties of an ABS-CG type algorithm 
for solving large linear systems},  International Conference on Computational
Mathematics,  Bangkock

Li C. (1999), {sl On a new algorithm of ABS-CG type for solving nonlinear
systems}, Numerical Methods and Computer Applications 20

Li S. (1997), {\sl The ABS methods for rank one perturbed linear system
of equation}, ABSC2, 45-52

Li S. (1995), {\sl On a variational characterization of the ABS algorithm},
preprint.

Li S. (1995), {\sl On the use of an ABS method in Karmarkar algorithm},
preprint.

Li S. and Zhang L. (1991), {\sl Some methods for base matrix Null(A) and 
numerical experiments}, ABSDAL.

Li Z.F. (1996), {\sl Restarted and truncated versions of ABS methods for
large linear systems: a basic framework}, QDMSIA 8/96.

Li Z.F. (1996), {\sl Modifying the Huang algorithm for parallel computation
of solutions to large and unstructured linear systems}, QDMSIA 25/96.

Li Z.F. (1996), {\sl Truncated block ABS methods for large linear systems},
QDMSIA 26/96.

Li Z.F. (1996), {\sl Applications of ABS methods to infinite linear systems},
QDMSIA 27/96.

Li Z.F. (1996), {\sl On the convergence of truncated ABS algorithms for large
linear systems}, QDMSIA 29/36.

Li Z.F. (1996), {\sl Towards ABS methods for ill-posed equality constrained
linear least squares problems}, QDMSIA 31/96.

Li Z.F., Deng N. and Wang L. (1996) {\sl ABS-type methods with ordering
strategies for solving nonlinear equations}, QDMSIA 12/96.

Longhino M. (1991), {\sl Un algoritmo ABS di tipo gradienti coniugati per 
minimi quadrati lineari}, Ms. Dissertation, University of Milan.

Longhino M., Spedicato E. and Vespucci M.T. (1992), {\sl On the
minimum residual conjugate gradient type ABS algorithm for
linear least squares}, QDMSIA 29/92.

Luksan L. and Spedicato E. (1998) {\sl Variable metric methods for
unconstrained optimization}, QDMSIA 7/98

Luo Y., Feng Z. and Liu J. (1997), {\sl On the conditions of the semilocal
convergence for nonlinear ABS algorithms}, ABSC2, 27-34

Meng F., Wang X. and He Z. (1997), {\sl A theorem on $W^m$ being $K_i$
admissible}, ABSC2, 53-56

Mirnia K. (1996), {\sl Iterative refinement in ABS methods}, QDMSIA 32/96.

Nicolai S. (1995), {\sl Metodi numerici per la soluzione delle equazioni
KT}, Ms. Dissertation, University of Bergamo.

Nicolai S. and Spedicato E. (1997), {\sl A bibliography of the ABS methods},
OMS 8, 171-183

Nicolai S. and Spedicato E. (1996), {\sl A bibliography on the ABS methods},
QDMSIA 3/96.

Oprandi M. (1987), {\sl Metodi  ABS per sistemi algebrici lineari},       
Ms. Dissertation, University of Bergamo.

Phua K.H. (1986), {\sl Solving sparse linear systems  by ABS methods}, 
Report 4/86, Computer Science Department, National University, Singapore.

Phua K.H. (1988), {\sl Solving sparse linear systems by an ABS method that 
corresponds to LU decomposition}, BIT 28, 709-718.

Ren Z. (1991), {\sl Solving optimization problems by combining the ABS methods 
and monotonic analysis}, ABSDAL.

Ronto N.I., Tuzson A. (1994) {\sl Construction of periodic solutions of 
differential equations with impulse effect}, Publ. Math. Debrecen, 44, 335-357.

Shahmorad S. (1994), {\sl Solution of linear least squares via the ABS 
algorithm}, Ms. Dissertation, Tabriz University.

Shan R. and Liu W. (1997), {\sl A finite terminating conjugate direction
method for linearly constrained optimization}, ABSC2, 35-44

Shi G. (1992), {\sl An ABS algorithm for generating nonnegative
solutions to linear systems}, ABSC1, 54-57. Also in ABSDAL.

Shi G. and Liu X. (1992), {\sl An ABS algorithm for the concave-minimum 
problem over  a  convex polyhedroid}, ABSC1, 58-64. Also in ABSDAL.

Spedicato E. (1983), {\sl On the symmetric algorithm in the ABS class of 
direct methods for linear systems}, QDMSIA 13/83.

Spedicato E. (1985), {\sl On  the  solution of linear  least  squares through 
the ABS class for linear systems}, Atti  Giornate AIRO, September 1985, 
89-98.

Spedicato E. (1986), {\sl Propriet\'{a}  di convergenza  degli algoritmi ABS 
per sistemi algebrici nonlineari}, Ricerca Operativa e Informatica, 
Milano, Franco Angeli, 57-65 (Bielli M. editor). 

Spedicato E. (1987a), {\sl A subclass  of variationally derived algorithms  
in the ABS class for linear systems},  Report 202, Mathematische 
Institute, University of W\"{u}rzburg.

Spedicato E. (1987b), {\sl Optimal conditioning parameter selection in the ABS 
class for linear systems}, Report 203, Mathematische Institute, University of 
W\"{u}rzburg.

Spedicato E. (1987c), {\sl ABS algorithms for sequential and parallel solution 
of linear and nonlinear algebraic equations}, QDMSIA 12/87.

Spedicato E. (1987d), {\sl Variationally derived algorithms in the ABS  class 
for linear systems}, Calcolo, 24, 241-246.

Spedicato E. (1988), {\sl Optimal conditioning parameter selection in the ABS 
class for linear algebraic systems}, Matem\'{a}tica Aplicada e Computacional, 
7, 187-199.

Spedicato E. (1989), {\sl A class of algorithms for sequential and parallel 
solution of algebraic linear and nonlinear systems}, 
Advances in Chinese Mathematics 18, 55-61 (in Chinese).
Also in Proceedings of the 7th 
Mathematical  Programming Symposium, Nagoya, 1986 (in English).  

Spedicato E. (1991a), {\sl ABS algorithms for general linear systems}, 
in Computer Algorithms  for Solving  Linear Algebraic Equations:  the State of 
the Art,  NATO  ASI  Series,  Vol. F77, 93-110, Springer-Verlag 
(Spedicato E. editor).

Spedicato E. (1991c), {\sl ABS algorithms for general linear 
systems}, QDMSIA 8/91.

Spedicato E. (1992a), {\sl A  class of sparse symmetric quasi-Newton updates}, 
Ricerca Operativa 22, 63-70, 1992. Also as QDMSIA 1/92. 

Spedicato E. (1992b), {\sl Special Issue: The First International
Conference on ABS Methods}, OMS 1, 169-281 (editor).

Spedicato E. (1992c), {\sl Proceedings of the  First International
Conference on ABS algorithms}, University of Bergamo (editor).

Spedicato E. (1992d), {\sl A review of ABS algorithms for general
linear systems}, OMS 1, 171-182.

Spedicato E. (1993), {\sl Ten  years of ABS  methods:  a review of
theoretical results and computational achievements}, Surveys 
on Mathematics for Industry, 3, 217-232.

Spedicato E. (1995a), {\sl ABS algorithms for linear systems and linear
least squares: theoretical results and computational performance},
Scientia Iranica 1, 289-303, 1995. Also as QDMSIA 2/94. 

Spedicato E. (1995b), {\sl ABS algorithms from Luoyang to 
Beijing}, QDMSIA 12/95.

Spedicato E. (1996), {\sl New ABS linear solvers and applications to
optimization}, PhD Dissertation, Dalian University of Technology.

Spedicato. (1997), {\sl ABS algorithms from Luoyang to Beijing}, OMS 8,
87-97

Spedicato E. (1998), {\sl De aliqua methodi Huangianae proprietate}, Vox
Latina 347-351, also as QDMSIA 4/91.

Spedicato E.  (1998), {\sl ABS algorithms for linear equations and linear
least squares: a review}, QDMSIA 15/98.

Spedicato E. (1999), {\sl ABS algorithms for linear equations and optimization:
a review of main results}, Proceedings  ISC99, Beirut, March 1999 (I. Moghrabi
and S. Kabbani editors), LAU University, 11-21
		  
Spedicato E. and Abaffy J. (1982), {\sl A class of direct methods for linear 
systems I: basic properties}, Report 82/4, IAMI, Milan.

Spedicato E. and Abaffy J. (1987) {\sl On the use of the ABS algorithms
for some linear programming problems}, Preprint, University of Bergamo.

Spedicato E. and Bodon E. (1987a), {\sl Numerical  experiments in the ABS 
class for nonlinear systems of algebraic equations}, Report SOFMAT 1/87, 
IAC, Rome.

Spedicato E. and Bodon E. (1987b), {\sl Application  of the ABS algorithm to 
linear least squares}, QDMSIA 13/87.

Spedicato E. and Bodon E. (1989a), {\sl Solving linear least squares by 
orthogonal factorization  and pseudoinverse computation via the modified 
Huang algorithm in the ABS class}, Computing, 42, 195-205. Also as QDMSIA
9/88.

Spedicato E. and Bodon E. (1989b), {\sl Biconjugate algorithms in the ABS 
class II: numerical evaluation}, QDMSIA  4/89.

Spedicato E. and Bodon E. (1990), {\sl Numerical  experiments of the implicit 
QR algorithm  of the ABS class  for nonlinear least squares}, Preprint, 
University of Bergamo.

Spedicato E. and Bodon E. (1991a), {\sl Numerical performance of conjugate 
gradient type ABS methods on ill-conditioned problems}, QDMSIA 13/91.

Spedicato E. and Bodon E. (1991b), {\sl Numerical  behaviour of  the implicit 
QR algorithm in the ABS class for nonlinear systems}, QDMSIA 16/91.

Spedicato E. and Bodon E. (1992), {\sl Numerical behaviour of the
implicit QR algorithm in the ABS class for linear least
squares}, Ricerca Operativa 22, 43-55. Also as QDMSIA 13/89.

Spedicato E. and Bodon E. (1993), {\sl Solution of linear least
squares via the  ABS  algorithm}, Mathematical Programming
58, 111-136. Also as QDMSIA 11/89.

Spedicato E., Chen Z. and Bodon E. (1996), {\sl ABS methods for KT
equations}, in Nonlinear Optimization and Applications, Plenum Press
(Di Pillo G. and Giannessi F. editors). Also as QDMSIA 13/95.

Spedicato E., Chen Z. and Deng N. (1993), {\sl The truncated difference
ABS-type algorithm for nonlinear systems of equations}, QDMSIA 21/93.

Spedicato E., Chen Z. and Deng N. (1994), {\sl A class  of difference ABS-type 
algorithms for a nonlinear system of equations}, 
Numerical Linear Algebra with Applications 13, 313-329. Also
as QDMSIA 6/91.

Spedicato E., Chen Z. and Deng N. (1996), {\sl A new nonlinear ABS-type
algorithm and its convergence property}, QDMSIA 33/96.

Spedicato E., Gao Z. and Yu H. (1990), {\sl The scaled ABS versions of some 
conjugate gradient methods for nonsymmetric linear systems}, QDMSIA 19/90.

Spedicato E. and Huang Z. (1992), {\sl An orthogonally scaled ABS
method for nonlinear least squares problems}, OMS 1, 233-242. Also as QDMSIA 
4/92.

Spedicato E. and Huang Z. (1995), {\sl Optimally stable ABS methods
for nonlinear underdetermined systems},  OMS 5, 17-26.
Also as QDMSIA 20/92.

Spedicato E. and Huang Z. (1996), {\sl Numerical experience with Newton-like
methods for nonlinear algebraic equations}, QDMSIA 96/15.

Spedicato E. and Li Z. (1997), {\sl A class of sparse symmetric
Quasi-Newton updates with least-change column-correction}, QDMSIA 97/17.

Spedicato E. and Li Z. (1998), {\sl Restarted and truncated versions of
ABS methods for large linear systems: a basic framework}, in Proceedings
of the First Workshop on Large Scale Scientific Computations, Varna, June
7-11, 1997

Spedicato E. and Li Z. (1998), {\sl Towards ABS methods for large-scale
sparse linear systems}, QDMSIA 98/1.

Spedicato E., Li Z. and Luksan L. (1997), {\sl Proceedings of the second
international conference on ABS algorithms. Beijing, June 1995}, QDMSIA 97/19.

Spedicato E. and Tuma M. (1993), {\sl Solving sparse unsymmetric
linear systems by implicit Gauss algorithm: stability}, QDMSIA 4/93.

Spedicato E. and Vespucci M.T. (1987), {\sl Computational performance of the 
implicit Gram-Schmidt  (Huang) algorithm for linear algebraic systems}, 
QDMSIA 4/87.

Spedicato E. and Vespucci M.T. (1992), {\sl A bibliography on the ABS
methods}, QDMSIA 19/92.

Spedicato E.  and  Vespucci M.T. (1993), {\sl Variations on the Gram-Schmidt 
and the Huang algorithms for linear systems: a numerical study}, 
Applications of Mathematics 2, 81-100. Also as QDMSIA 21/89.

Spedicato E. and Xia Z. (1990a), {\sl Application of ABS algorithms to 
constrained optimization I:  linear equality constraints}, QDMSIA 7/90.

Spedicato E.  and  Xia Z.  (1990b), {\sl An  approach  for polyhedral 
approximations to subdifferentials of finite convex functions}, QDMSIA 17/90.

Spedicato E. and Xia Z. (1992a), {\sl Finding general solutions of the
Quasi-Newton equation via the ABS approach}, OMS 1, 243-252.
Also as QDMSIA 10/91.

Spedicato E. and Xia Z. (1992b), {\sl A class of descent algorithms for
nonlinear function minimization with linear equality
constraints}, OMS 1, 265-272. 

Spedicato E. and Xia Z. (1992c), {\sl On  some  ABS methods for the 
computation of the pseudoinverse}, Ricerca Operativa 22, 
35-41. Also as QDMSIA 24/90.

Spedicato E. and Xia Z. (1994), {\sl ABS algorithms for nonlinear optimization}, 
in Algorithms for Continuous Optimization: the State of the Art, Kluwer, 
333-356 (Spedicato E. editor).

Spedicato E. and Xia Z. (1998), {\sl ABS algorithms for optimization: a 
review}, QDMSIA 16/98.

Spedicato E. and Xia Z. (1999), {\sl ABS formulation of Fletcher implicit LU
method for factorizing a matrix and solving linear equations and  LP problems},
QDMSIA 1/99

Spedicato E. and Yang Z. (1990), {\sl Optimally conditioned scaled ABS algorithms 
for linear systems},  JOTA 67, 141-150. Also as QDMSIA 15/89

Spedicato E. and Zhao J. (1992), {\sl On sparse symmetric and positive 
definite quasi-Newton updates}, QDMSIA 2/92.

Spedicato E. and Zhao J. (1993), {\sl Explicit general solution of the
Quasi-Newton equation with sparsity and symmetry}, OMS 2, 311-319.
Also as QDMSIA 11/93.

Spedicato E. and Zhu M. (1994), {\sl On the generalized implicit
LU algorithm of the ABS class for linear systems}, QDMSIA 3/94.

Spedicato E. and Zhu M. (1996), {\sl Reformulation of the ABS algorithm
via full rank Abaffians}, QDMSIA 4/96.

Spedicato E. and Zhu M. (1996), {\sl Efficient solution of  linear systems
after rank-one correction via the ABS algorithm with full rank Abaffians},
QDMSIA 5/96.

Spedicato E. and Zhu M. (1999), {\sl A generalization of the implicit LU
algorithm to an arbitrary initial matrix}, Numerical Algorithms 20, 343-351

Spedicato E., Esmaeili H. and Xia Z. (1999), {\sl A review of ABS
algorithms for linear real and Diophantine equations and optimization},
Proceedings  XVI CEDYA and VI CMA Conferences, Las Palmas de Gran Canaria,
September 1999 (Montenegro R.,  Montero G. and Winter G. editors), CEANI,
11-23

Spedicato E., Xia Z. and Gao Z. (1990), {\sl The basic block ABS method I: 
formulation and theoretical properties}, QDMSIA 15/90.

Spedicato E., Xia Z. and Zhang L. (1997), {\sl The implicit LX method of
the ABS class}, OMS 8, 99-110, also as QDMSIA 2/96.

Spedicato E., Xia Z., Zhang L. and Mirnia K. (1996), {\sl ABS algorithms
for linear equations and applications to optimization}, QDMSIA 30/96.

Spedicato E., Xia Z. and Zhang L. (1999), {ABS algorithms for solving
linearly constrained optimization problems via the active set strategy},
QDMSIA 16/99

Spedicato E., Zhang L. and Xia Z. (1999), {\sl Solving linear matrix
equations via the ABS method and applications to Quasi-Newton equations},
QDMSIA 12/99

Tan Z. (1991), {\sl A new algorithm with finite convergence to solve a class 
of quadratic programming constructed by the ABS methods}, ABSDAL.

Tan Z. and Chen Y. (1991), {\sl An improvement on suboptimization methods for 
convex programming constraints}, ABSDAL.

Tan Z. and Gao Y. (1991), {\sl A quasi-gradient projection method for 
nonlinear constrained systems}, ABSDAL.

Tan Z. and Gao Y. (1992), {\sl An improved gradient projection method for  
nonlinear constrained systems}, ABSC1, 65-77.

Tuma M. (1990), {\sl Utilization of matrix sparsity in ABS methods for solving 
systems of linear equations},  QDMSIA 21/90.

Tuma M. (1993), {\sl Solving sparse unsymmetric sets of linear equations
based on implicit Gauss projections}, Report 556, Institute of Computer
Science, Academy of Sciences, Prague.

Tuma M. (1996), {\sl On orthogonalization procedures in some iterative algorithms
for systems of linear equations}, QDMSIA 8/96.

Vespucci M.T. (1990), {\sl  On the implicit Gauss-Cholesky algorithm of the 
ABS class}, Contributed Papers, NATO ASI on Computer Algorithms for  
Solving Linear Algebraic Equations: the State of the Art, Il Ciocco, 1990, 
Cooperativa Studium Bergomense, 186-205 (Spedicato E. editor).

Vespucci M.T. (1991a), {\sl The use of the ABS algorithms in truncated Newton
methods for linear optimization}, Colloquia Mathematica Societatis J\'{a}nos 
Bolyai 59,  Numerical Methods, 81-91 (Stoyan G. editor).

Vespucci M.T. (1991b), {\sl An ABS implicit Gauss-Cholesky code for ND matrices}, 
Manual 1/68, Progetto Finalizzato Sistemi Informatici e Calcolo Parallelo, 
IAC, Roma.

Vespucci M.T. (1991c), {\sl Truncated Newton methods based on the ABS class}, 
Ph.D. Dissertation, Hatfield Polytechnic. 

Vespucci M.T. (1992), {\sl The ABS preconditioned conjugate gradient algorithm}, 
QDMSIA 6/92. Also in ABSC1, 78-94.

Vespucci M.T., Yang Z., Feng E., Yu H. and Wang X. (1992), {\sl A bibliography
on the ABS methods}, OMS 1, 273-281.

Vitali A. (1989), {\sl Algoritmi vettoriali e paralleli della classe ABS 
per la soluzione  di sistemi lineari}, Ms. Dissertation , University of Milan.

Xia Z. (1990a), {\sl Application  of  ABS algorithms to constrained 
optimization II: linear inequality constraints}, QDMSIA 29/90.

Xia Z. (1990b), {\sl An algorithm  for minimizing a class of quasidifferentiable 
functions}, QDMSIA 12/90.

Xia Z. (1991a), {\sl Quadratic programming  via the ABS algorithm I: equality 
constraints}, QDMSIA 1/91.

Xia Z. (1991b), {\sl The  ABS  class and quadratic programming II: general  
constraints},  QDMSIA 5/91.

Xia Z. (1992a), {\sl Finding  subgradients  or  descent directions of convex 
functions by external polyhedral approximation of subdifferentials}, OMS 
253-264.

Xia Z. (1992b), {\sl Derivation of some linear transformations  via
the ABS algorithm}, ABSC1, 95-110.

Xia  Z. (1994), {\sl An  efficient ABS algorithm for linearly
constrained optimization}, QDMSIA 15/94.

Xia Z. (1995a), {\sl ABS  formulation and generalization of the
reduced gradient method}, QDMSIA 3/95.

Xia Z. (1995b), {\sl ABS formulation  and generalization of the
interior point method}, ABSC2.

Xia Z. (1995c), {\sl ABS reformulation of some versions of the simplex
method for linear programming},  QDMSIA 10/95.

Xia Z. (1998), {\sl Further ABS formulation of reduced gradient methods},
QDMSIA 2/98.

Xia Z. and Gao Z. (1990), {\sl  The  basic block ABS methods II: application 
to the Rutherford matrix equation}, QDMSIA 16/90.

Xia Z. and Liu Y. (1993), {\sl On ABS descent directions for constrained
minimization with linear equality constraints}, Journal of Dalian
University of Technology 35, 255-263.

Xia Z. and Zhang L. (1999), {\sl A modified outer approximation method for
solving concave programming problems via the ABS algorithm}, QDMSIA 21/99

Xia  Z.,  Liu  Y. and Zhang L. (1992), {\sl Application of a
representation of ABS updating  matrices  to  linearly
constrained optimization I}, Northeast Operational Research 7, 1-9.

Xia Z., Zhang L. and Liu Y. (1995) {\sl An ABS algorithm for minimizing a
nonlinear function subject to linear inequalities}, ABSC2.

Xiu N. and Chen G. (1991), {\sl Optimality of the ABS matrices 
update formula}, Preprint, Qufu Normal University.

Xue Y. (1997), {\sl A practical method for solving linear and nonlinear
least squares problems}, ABSC2, 17-26

Yang Z. (1989a), {\sl ABS algorithms  for solving  certain systems 
of indefinite equations}, QDMSIA 6/89.

Yang Z. (1989b), {\sl The properties of iterative point in ABS algorithms 
and the related generalized inverse matrix}, QDMSIA 7/89. 

Yang Z. (1989c), {\sl On the numerical stability of the Huang and the 
modified Huang algorithms and related topics}, QDMSIA 8/89.

Yang Z. (1989d), {\sl An  ABS subclass with optimal conditioning parameter 
selection}, QDMSIA  10/89.

Yang Z. and Zhu M. (1994), {\sl The practical ABS algorithm for large
scale nested dissection  linear system}, QDMSIA 11/94.

Wang X., Cai J. and Li X. (1995), {\sl On the application of the ABS
algorithm in the simplex method}, ABSC2.

Wang X., Feng E. and  Liu  X. (1995), {\sl ABS reformulation of
Lemke's method for the linear complementarity problem}, ABSC2.

Wang Z. (1997), {\sl ABS formulation of Bunch-Parlett decomposition},
ABSC2, 9-16

Zeng K. and Jiang G. (1991), {\sl An application of the ABS algorithms to 
solve nonlinear systems of equations}, Preprint, South-Western Jiaotong 
University.       

Zeng L. (1989), {\sl An  ABS  algorithm for sparse linear systems}, 
QDMSIA 9/89.

Zeng L. (1991), {\sl The parallel scaled ABS class for sparse systems of 
linear equations}, ABSDAL.

Zeng L. (1997), {\sl A parallel ABS method for sparse nonlinear
algebraic equations}, ABSC2, 1-8

Zhang J. and Zhu M. (1995a), {\sl A symmetric ABS algorithm for some matrix
matrix equation and its minimum properties}, Report MA95-19, City
University, Hong Kong.

Zhang J. and Zhu M. (1995b),  {\sl Solving a least change problem under the weak
secant condition by an ABS approach}, Report MA95-20, City University,
Hong Kong

Zhang J. and Zhu M. (1997), {\sl Least-change properties of ABS methods and their
application in secant-type updates}, OMS 8, 111-132

Zhang L. (1992a), {\sl Application of the ABS algorithm  to Quasi-Newton 
methods for linearly constrained optimization problems}, ABSC1 119-130.
Also in ABSDAL.

Zhang L. (1992b), {\sl A method for finding a  feasible point  
of inequalities}, ABSC1 131-137. Also in ABSDAL.

Zhang L. (1992c), {\sl A class of quasi-Newton formulae with $ZZ^T$
factorization}, Proceedings of the First Youth Academic Meeting of Liaoning, 
1992, 220-224, Northeast University Press.

Zhang L. (1995a), {\sl An algorithm for the least Euclidean norm
solution of a linear system of inequalities via the  Huang
ABS algorithm and the Goldfarb-Idnani strategy}, QDMSIA 2/95.

Zhang  L. (1995b), {\sl Relations between projection  matrices
and ABS updating matrices},  QDMSIA 5/95.

Zhang L. (1997), {\sl On the ABS  algorithm with singular initial
matrix and its application to linear programming}, OMS 8 143-156
(also as QDMSIA 8/95).

Zhang L. (1995d), {\sl Updating of Abaffian matrices under perturbation
in W and A}, QDMSIA 15/95.

Zhang L. (1999), {\sl Computing inertias of KKT matrix and reduced Hessian
via the ABS algorithm for applications to quadratic programming}, 
QDMSIA 7/99

Zhang L., Feng E. and Jiang Z. (1992), {\sl A class of projection 
methods for linearly constrained optimization problems}, ABSC1, 138-149.
Also in ABSDAL.

Zhang L., Feng E. and Xia Z. (1993), {\sl A theorem on the steepest descent
direction for linearly constrained problem,} Journal of Mathematical Research 
and Exposition 13, 358.

Zhang L. and Xia Z. (1993), {\sl A modified version of ABS algorithms}, 
Journal of Mathematical Research and Exposition 13, 26.

Zhang L. and Xia Z. (1994), {\sl Comparing ABS algorithm with Fletcher
and GI methods}, Northeast Operations Research, 8, 36-38.

Zhang L. and Xia  Z.H. (1995), {\sl Application of the implicit  LX
algorithm to the simplex method}, QDMSIA 9/95.

Zhang L., Xia Z. and Feng E. (1998), {\sl Introduction to ABS Methods in
Optimization}, Monograph, Dalian University of Technology Press
							       
Zhang L., Xia Z. and Wang X. (1994), {\sl A version of updating ABS matrices
in optimization}, Northeast Operations Research, 8, 51-54.
		    
Zhang J. and Zhu M. (1997), {\sl Least change properties of ABS methods and their
application in secant-type updates}, OMS 8, 111-132

Zhao J. (1981), {\sl Huang's method for solution of consistent linear 
equations and its generalization},  Journal of Numerical Mathematics of 
Chinese Universities, 3, 8-17.

Zhao J. (1989), {\sl A class of direct methods for solving linear 
inequalities}, Journal of Numerical Mathematics of Chinese Universities
3, 231-38.

Zhao J. (1991), {\sl ABS algorithms  for solving linear inequalities}, 
QDMSIA 21/91.

Zhao J. (1995), {\sl A note on the ABS algorithm for solving linear systems},
Journal of Nanjing University 12, 32-36.

Zheng H., Fei P. and Fan R. (1991), {\sl The parallel Huang algorithm of the 
ABS class on MIMD parallel system}, Preprint, University of Wuhan.

Zhu M. (1987), {\sl The implicit  $LL^T$ algorithm for sparse nested 
dissection linear systems}, Technical Report 196, Numerical 
Optimisation Centre, Hatfield Polytechnic.

Zhu M. (1989a),  {\sl The Q-quadratic  convergence for a subclass of 
nonlinear ABS algorithms}, QDMSIA 22/89.

Zhu M. (1989b), {\sl A generalization of nonlinear ABS algorithms and its 
convergence properties}, Journal of Beijing Polytechnic University, 
15, 21-26. 

Zhu M. (1990), {\sl On the ABS algorithms for perturbed linear systems}, 
Journal of Computer Mathematics 29, 181-194. Also as
Report 198, Numerical Optimisation Centre, Hatfield Polytechnic, 1987. 

Zhu M. (1991a), {\sl Finding  general  solution to perturbed linear systems 
and its application in optimization}, Communication on 
Applied and Computational Mathematics, 5, 61-67. Also as QDMSIA 5/90

Zhu M. (1991b), {\sl A note on the condition for convergence  of
nonlinear ABS method}, Calcolo 28, 307-314. Also as QDMSIA 17/91. 

Zhu M. and Zhang J. (1995), {\sl Solving a least change problem
under the weak secant condition by the ABS approach}, ABSC2.
\end{document}